\documentclass[10pt,journal,cspaper,compsoc]{IEEEtran}
\ifCLASSOPTIONcompsoc
 \else
\fi
\usepackage{algorithm2e}

\ifCLASSINFOpdf
  \usepackage[pdftex]{graphicx}
  
\else
  \usepackage[dvips]{graphicx}
\fi
\usepackage{amssymb}
\usepackage[cmex10]{amsmath}
 \usepackage{epstopdf}
\usepackage{picins}
\newtheorem{theorem}{Theorem}

\newtheorem{definition}[theorem]{Definition}

\newtheorem{proposition}[theorem]{Proposition}

\newtheorem{remark}[theorem]{Remark}

\newcommand{\real}{\ensuremath{\mathbb{R}}}

\newcommand{\innerd}[2]{\left\langle\!\left\langle #1,#2 \right\rangle\!\right\rangle}

\usepackage{hyperref}

\usepackage[colorinlistoftodos]{todonotes}

\begin{document}
%
\title{Gauge Invariant Framework for Shape Analysis of Surfaces}
%
%
%
%

\author{Alice~Barbara~Tumpach,
        Hassen~Drira,
        Mohamed~Daoudi,~\IEEEmembership{Senior,~IEEE}, and~Anuj~Srivastava,~\IEEEmembership{Senior,~IEEE}
\IEEEcompsocitemizethanks{\IEEEcompsocthanksitem A.B. Tumpach is with Laboratoire Paul Painlev\'e (U.M.R. CNRS 8524) University of Lille1, France.   \protect\\
E-mail: Barbara.Tumpach@math.univ-lille1.fr
\IEEEcompsocthanksitem H. Drira is with Institut Mines-Telecom/Telecom Lille CRIStAL (UMR CNRS 9189).
E-mail: hassen.drira@telecom-lille.fr
\IEEEcompsocthanksitem M. Daoudi is with Institut Mines-Telecom/Telecom Lille CRIStAL (UMR CNRS 9189).
E-mail: mohamed.daoudi@telecom-lille.fr
\IEEEcompsocthanksitem A. Srivastava is with the Department of Statistics, Florida State University, Tallahassee, USA.
E-mail: anuj@fsu.edu

}
\thanks{}}

%
%

\markboth{IEEE PATTERN ANALYSIS AND MACHINE INTELLIGENCE}%
{Tumpach \MakeLowercase{\textit{et al.}}: Gauge Invariant Framework for Shape Analysis of Surfaces}

\IEEEcompsoctitleabstractindextext{%
\begin{abstract}
This paper describes a novel framework for computing geodesic 
paths in shape spaces of spherical surfaces under an elastic Riemannian metric. 
The novelty lies in defining this Riemannian metric directly on the quotient (shape) space, 
rather than inheriting it from pre-shape space, and using it to formulate  a
path energy that measures only the normal components of velocities along the path.
In other words,
this paper defines and solves for geodesics directly on the shape space and avoids complications
resulting from the quotient operation. 
This comprehensive framework is invariant to arbitrary parameterizations of surfaces along 
paths, a phenomenon termed as gauge invariance. 
Additionally, this paper makes a link between different elastic metrics used in the computer science literature on one hand, 
and the mathematical literature on the other hand, and provides a geometrical interpretation of the terms involved.
Examples using real and simulated 3D objects
are provided to help illustrate the main ideas. 
\end{abstract}

\begin{keywords}
3D surfaces, Riemannian metric, geodesics.
\end{keywords}}

\maketitle


\IEEEdisplaynotcompsoctitleabstractindextext

%
\IEEEpeerreviewmaketitle

        	\begin{figure*}[!ht]
 		\centering
 		\includegraphics[width=16cm]{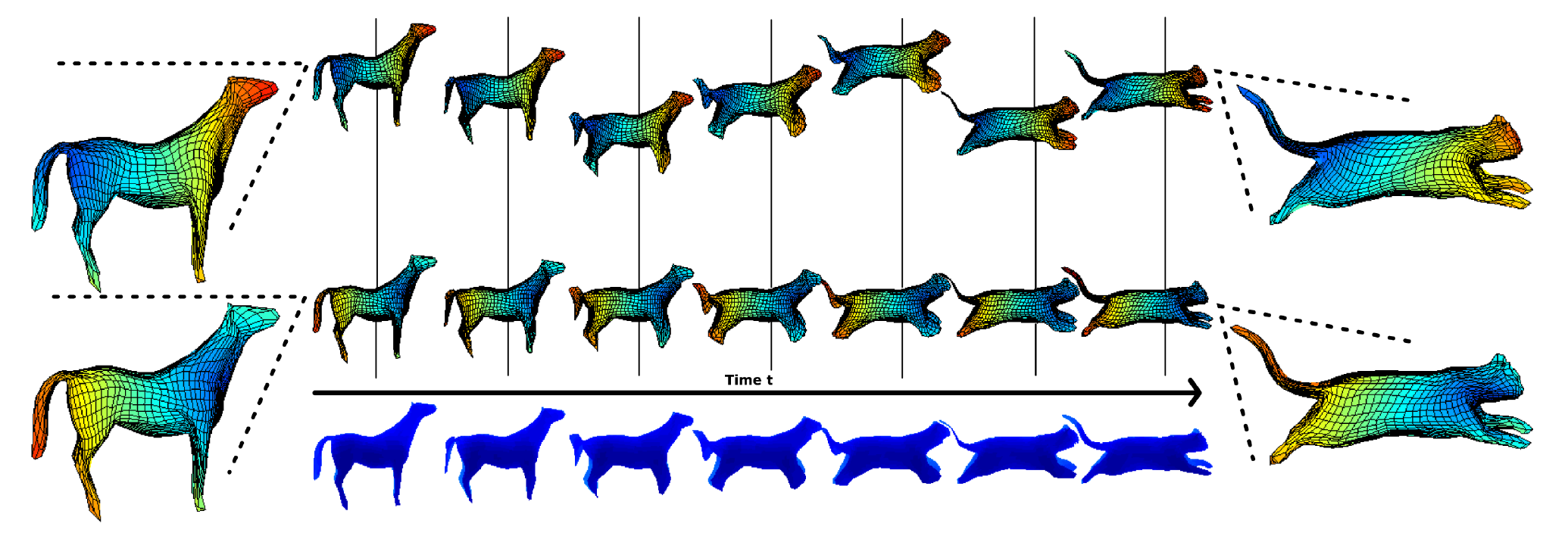}
 		\caption{\scriptsize
        Two paths in $\mathcal{F}$ with the same sequence of shapes but with different reparameterizations between the corresponding shapes.}
		\label{fig_illustration}
		\end{figure*}

\section{Introduction}
In this paper we seek a framework for analysing shapes of a certain class of 
3D objects. Although the general goal in shape analysis is to develop tools for full statistical analysis -- statistical averaging, 
finding principal modes of variations in a population, and shape  classification, we restrict to more basic
goals of quantifying shape differences and generating deformations. 
While there have been many efforts in shape analysis of 3D objects, the 
problem is far from solved and the current solutions face many technical and 
practical issues. For instance, many general techniques for shape analysis rely on quantifying shape 
differences by spatially matching geometric features
across objects. Therefore,  it becomes important to establish a correspondence of parts between objects, i.e. which part in 
one object corresponds to which part in the other? This was an 
important bottleneck in a majority of previous efforts on 3D shape analysis where the correspondence (or registration) of objects was 
either presumed or solved as an independent pre-processing step. More recently, there
has been progress in establishing frameworks that formulate the registration and comparison problems 
jointly.
These newer frameworks, using techniques from differential geometry,  focus on shape analysis of parameterized surfaces 
and treat the problem of shape comparison as the problem of computing geodesic paths in shape spaces
under a chosen metric.  Here shapes are compared using a Riemannian metric on a \textit{pre-shape space} $\mathcal{F}$ consisting 
of embeddings or immersions of a model manifold (like the sphere, or the disc) into 
the 3D Euclidean space $\mathbb{R}^3$.  Two embeddings correspond to the same shape in $\mathbb{R}^3$ if and only if 
they differ by an element of  a shape-preserving transformation group, such as rigid motion, scaling, and reparameterization.
The shape space is therefore the quotient space of the pre-shape space by these shape-preserving groups. If the Riemannian metric on the 
pre-shape space is preserved by the action of the shape-preserving group then it induces a Riemannian metric on the quotient space. 
The construction of geodesics in shape space provide optimal deformations between surfaces and is a very important tool 
in {\it statistical} analysis of shapes. 
Interestingly, the problem of registration is handled using parameterizations of surfaces such that 
the points denoting the same parameter values on two objects are considered registered. \\

While these geometric 
ideas are powerful and comprehensive, there are two important issues that one needs to deal with: 
(1) the choice of Riemannian metric to define geodesics, geodesic lengths, and the eventual shape metric,  and (2)  the task
of computing geodesic paths between arbitrary shapes. In terms of the first issue, the choice of a metric, an important requirement 
is that the metric should be invariant to action of the reparameterization group, to enable a well-defined distance on the 
eventual quotient space or the shape space of surfaces.  There is a related requirement for the shape analysis to be invariant to 
parameterizations of objects since parameterizations are only artificial impositions designed to help navigate along objects. 
The physical intuition we have is that shape tools, such as the deformation (path or geodesic) from one shape to another,
are physical processes that are independent of the way surfaces may be parameterized. 
These dual requirements rule out the use of commonly-used quantities such as the $L^2$ norm on the space ${\cal F}$ directly. 
In terms of the second issue, the lack of standard metrics makes it complicated to compute geodesic paths even when 
the underlying manifold is a vector space, and one needs numerical algorithms for approximating geodesic paths. Next, we 
present a summary of the past work on these two issues and outline motivations for the current paper.

\subsection{Motivation and Past Work}

Our goal is to develop tools for analyzing shapes of two-dimensional
surfaces with certain local constraints (smoothness, no-holes, etc).
The main difficulty in comparing shapes of such surfaces is that
there is no preferred parameterization that can be used for
registering and comparing features across surfaces. Since the shape
of a surface is invariant to its parameterization, one would like an
approach that yields the same result irrespective of the
parameterization. 

Furthermore, we are not only interested in the comparison and matching of two shapes, 
but also in the deformation processes that may transform one shape into  another, i.e. metamorphosis. 
To be physically meaningful, the evolution from one shape to another
should be independent of the way surfaces may be parameterized. 
Our approach 
to shape analysis presented in this paper was therefore initiated by the following question~:
\textit{What is the natural framework where one can measure deformations of shapes independently 
of the way shapes are parameterized?} 
As a motivating example, the sequence of shapes displayed in Fig. \ref{fig_illustration} (bottom) 
denotes a path where a horse is transformed into a jumping cat. 
During the transformation process, only the change of shape, drawn in the bottom line as a sequence of blue surfaces, is relevant to us. 
How the surfaces may be parameterized 
during the metamorphosis has no importance in our context. 
To emphasize this idea, two paths of parameterized surfaces corresponding to 
the same  transformation process are displayed in the top two rows.  We would like  a framework where the physical quantities 
measured on the path of shapes, 
such as its length or its energy, are independent of the parameterizations of surfaces 
along the transformation process. In particular, in Fig. \ref{fig_illustration}, 
the two paths of parameterized surfaces
corresponding to the same transformation process should have the same length. 
Note that the surfaces along the second path are obtained by applying a \textit{different reparameterization
at each time step}  to the surfaces along the first path. 

Let us emphasize that we are not only  interested in how \textit{far} the horse and the jumping cat are from each other, in other words in a quantity like a distance measuring
the minimal cost needed to deform the horse into a cat. But, given a metamorphosis between these two shapes, we are also interested in measuring 
its \textit{length} on one hand, and its \textit{energy} on the other hand, independently of the parameterizations  of the transformation process that may have been used to create this metamorphosis. Recall that the length of a path is the integral of the norm velocity function with respect to time and has the dimension of a distance. The energy is the integral of the square of the norm velocity function with respect to time, hence has the dimension of the square of a distance divided by time.

Let us now summarize past work on related subjects. The initial set of papers developed algorithms for geodesic deformations between surfaces while using the given registration of points. 
They  compute geodesics between shapes, under isometric deformations, while assuming the registration (or parameterization) as given. 
Windheuser et al. \cite{Windheuser}  proposed to find a geometrically consistent matching of 3D shapes which minimizes an
elastic deformation energy but use a linear interpolation between registered pairs of points in $\mathbb{R}^3$ to compute geodesic paths.
Another paper by Kilian et al. \cite{kilian-2007-gmss}  represents parameterized surfaces by discrete triangulated meshes, assumes a Riemannian metric on the space of such meshes, and computes geodesic paths between given meshes. The main limitation 
here is that it assumes the correspondence between points across meshes. That is, we need to know beforehand which point on one mesh matches with which point on the second mesh. 
The same limitation holds for the paper by Heeren et al. \cite{Heeren-etal:2012} also. In contrast, we would like to remove the reparameterization variability so that different surfaces with the same shape but different parameterizations have zero distance between them.
 
Motivated by progress in shape analysis of curves \cite{Younes98, Srivastava11}, Kurtek et al. \cite{KurtekCVPR2010}, \cite{KurtekPAMI2012} introduced a new representation, termed a {\it $q$-map} of surfaces such that the 
$L^2$ distance in this representation space is invariant to simultaneous reparameterizations of surfaces. 
For convenience of the reader, we recall the definition of the $q$-map but we will not use it in the present paper.
Let $f: \mathbb{S}^2 \to \mathbb{R}^3$ denote a smooth parameterized surface and ${\cal F}$ be the set of such surfaces. 
Then, this $q$-map is given by $f \mapsto q$ where $q(s) = \sqrt{r(s)} f(s)$ and $r(s)$ is the area multiplication factor 
of $f$ at $s \in \mathbb{S}^2$. 
They defined 
a Riemannian metric on the space of parameterized surfaces by {\it pulling back} $L^2$ metric under the $q$-map, and 
used a {\it path-straightening} algorithm to compute geodesic paths between given surfaces in a pre-shape space. 
This path-straightening is an iterative algorithm that updates an arbitrary initial path using the gradient of the energy function
mentioned above,
until the path converges
to a geodesic. The energy gradient  is approximated numerically using an (approximate) finite basis for ${\cal F}$. 
To remove the effects of original parameterizations, and to obtain geodesics in the shape space, 
they solve for an optimal reparameterization of one of the surfaces, under the same energy.  
There are  several other papers, including \cite{litke-etal:2005}, that focus exclusively on the task of finding optimal correspondence between 3D objects, either using physically-motivated energies or Riemannian metrics. Due to the use of gradient-based searches, these methods and previously mentioned papers do not guarantee a global solution, either for geodesics or for registration. In path-staightening, however, it can be shown that a path that is a local minimum of the path energy is a geodesic path, albeit not the shortest geodesic. To our knowledge, very few methods guarantee a globally-optimal solution to the problem of finding geodesics in shapes spaces of surfaces.
Although \cite{KurtekCVPR2010}
was the first to provide a geometric framework for joint registration-comparison problem, the Riemannian metric used there
has a limitation that it was not translation invariant.

To handle the translation issue mentioned above, Jermyn et al. \cite{JermynECCV2012} introduced a comprehensive 
Riemannian metric that  has several improvements, including the fact that it was translation invariant and 
allows some physical interpretations in its use. This metric, given later in Eqn. \eqref{eqn:elastic-metric}, has terms that can be
interpreted as measurements of bending, stretching, and changes in local curvatures of surfaces. (We elaborate on 
this topic later in Section \ref{sec:geometric-interp}.)
It has been termed an {\it elastic metric} because it is invariant to reparameterizations  and
the physical interpretations associated with it. 
Although \cite{JermynECCV2012} introduced this metric, it did not use the full metric to compute geodesic paths. 
Instead, it defined a new map, termed the square-root normal field, given by $q(s) = \sqrt{r(s)} {n}_f(s)$ 
where ${n}_f(s)$ denotes the unit normal to the surface at the point $s \in \mathbb{S}^2$.
The square-root normal field has the property that the {\it last two terms} of the elastic metric transform to the $L^2$ metric under the map
$f \mapsto q$, for some weighting of last two terms in the metric. The first term of the metric is discarded in this analysis. 
The transformation to $L^2$ metric is useful since one can apply some common tools from 
Hilbert space analysis to this problem, including the optimization over the reparameterization 
group for optimal registration, but this mapping $f \mapsto q$ is not onto and, hence, not invertible. 
The optimization step is challenging  because  the reparameterization group is an infinite-dimensional 
Fr\'echet Lie group, and the exponential map is not a local diffeomorphism.
Since the first term of the elastic metric introduced in \cite{JermynECCV2012} is not used by Jermyn et al., it can result in zero shape distance between two surfaces
that actually have different shapes. For example, a thin-tall cylinder and a fat-short cylinder, with same surface 
areas and unit normals, will have zero shape difference under this framework.

Another line of work in shape analysis comes from Michor et al. \cite{MM}, Bauer et al. \cite{BHM3}, \cite{BHM5}, \cite{BHM4} 
and Fuchs et al. \cite{FJSY09} (see also \cite{BBM} for an overview of a lot of mathematical results in this area). 
Different types of metrics have been studied~: Sobolev metrics in \cite{BHM4}, curvature weighted metrics in \cite{BHM3}, almost local metrics in \cite{BHM5},
metrics mesuring the deformations of the interiors of shapes in \cite{FJSY09}.
Let us mention that the first two terms of the metric we use in the present paper fit in the general study laid out  in \cite{BHM4}, and are related to the metrics 
 studied in \cite{Ebin}, \cite{Freed}, \cite{GilMedrano}  (in a sense that we will make clear in Section~\ref{sec:geometric-interp}).
In this set of papers, the idea is to replace the problem of solving 
the geodesic equation on shape space by the equivalent problem of solving the equation for horizontal geodesics in the pre-shape space. 
A geodesic in pre-shape space is horizontal if it is orthogonal to the orbits of the reparameterization group. One task in this strategy is therefore to 
compute the horizontal space on which the quotient map is an isometry,
or equivalently solve a minimization problem for the infinitesimal energy.
 Depending on the Riemannian metric on the pre-shape space, this task may be computationally trivial or extremely difficult to implement 
 (for metrics used in \cite{BHM3} and \cite{BHM5} it is just the space of normal vector fields, but for metrics used in \cite{BHM4} and \cite{FJSY09}  it involves the inversion of a pseudo-differential operator). 
 Another main contribution of these authors is to give sufficient conditions under which the Riemannian metric induced on shape space separates points, i.e. gives a non-zero geodesic distance between pairs of different shapes (a condition that is necessary to make shape comparison). It is worth noting  that, in this infinite-dimensional context, vanishing geodesic distance is a common phenomenon (as was first highlighted in \cite{MM}).
 For the metric we use, non-vanishing geodesic distance is guaranteed by the non-vanishing geodesic distance on the space of Riemannian metrics proved in \cite{Clarke} (at least on pairs of shapes inducing different pull-back metrics on the sphere, which is what we are interested in practice).

To summarize, the past approaches involving Riemannian geometry have tended to perform shape analysis in two steps. 
First, they select a representation space, or a pre-shape space, for objects of interest  -- curves \cite{Younes98, Srivastava11,BBCMM} and surfaces
\cite{KurtekPAMI2012,JermynECCV2012,BHM3, BHM5,BHM4} -- and impose a Riemannian structure on it ensuring that the 
actions of shape-preserving groups are by isometries. Next, they inherit this metric to the quotient space of the 
pre-shape space modulo the requisite groups, called the shape space, and seek geodesics between objects in 
this shape space. The task of inheriting Riemannian metrics to quotient spaces is complicated because 
reparameterization groups are Fr\'echet  Lie groups and the process of inheriting a metric requires closed orbits, as can be 
seen in \cite{LRK, Younes98,BBCMM, BHM4}, etc. Even though endowing shape space with a Riemannian metric (with positive distance function) seems to be a good approach, 
inducing this metric by a Riemannian metric on pre-shape space leads to difficulties that one would like to avoid 
(recall that we are only interested in shapes and not in the way they are parameterized). We will pursue a different strategy where the Riemannian metric is 
directly imposed on the quotient space, thus avoiding the need to satisfy conditions for inheriting metrics from the 
pre-shape space or computing an abstract horizontal space. Motivated by an easy implementation of the metrics, we take the point of view 
where the space of interest is the space of normal vector fields (in contrast with the horizontal space of a Riemannian submersion). Let us emphasize that there is no restriction in doing so : any Riemannian metric on shape space can by expressed 
as a metric defined on normal vector fields.

\subsection{Goals and Contributions}
Now we present the goals and contributions of this paper, and start by revisiting the question:  
\textit{What should be a good Riemannian metric on shape space~?} A good Riemannian metric on shape space should be such that~:
(1)  it induces a positive distance function on shape space, i.e. the infimum of the lengths of paths connecting two different shapes should be non-zero~; 
(2) the distance between two shapes should be independent of the way the two shapes are parameterized~; and, 
(3) the length of a path of shapes should be independent of the way shapes along the path are parameterized.
The last point should be thought of as the natural generalization of the fact that, on a finite-dimensional Riemannian manifold, the length of a curve is independent of the way the curve is parameterized. 
It should be true for any path (not only for geodesics), and is called {\it gauge invariance}.
Indeed the use of parameterized surfaces in order to measure the deformation of a shape can be compared to the use of a gauge. 
Let us comment on Fig. \ref{fig_illustration} in order to illustrate this idea. 
Each column depicts an orbit under the reparameterization group for the corresponding surface,  
the surfaces in a given orbit correspond to the same shape but with different parameterizations. A path of shapes can be lifted in many ways to a path of parameterized 
surfaces. In Fig. \ref{fig_illustration} 
two lifts of the bottom line path are depicted. The first path connects parameterized surfaces with 
different ``heights" in the fibers. This is made to emphasize that the variations of the 
``height" (i.e. of the parameterization) in the fibers should not influence the value of the length of the path of shapes.

The main contributions of this paper are following: 
\begin{itemize}
\item The proposed method achieves gauge invariance, i.e. 
the lengths of paths (geodesics or otherwise) measured under this metric are invariant to arbitrary reparameterizations 
of shapes along these paths (in particular, the two paths in Fig.\ref{fig_illustration} have the same length). 
\item 
It uses an elastic metric that accounts for any deformation of patches to define and compute geodesic paths between 
given objects in the shape space, and it presents a geometric interpretation of the different terms involved in this metric.
\item 
By defining a metric directly in the shape space,  it avoids the optimization step over the reparameterization group and 
difficult mathematical issues arising from inheriting a metric from pre-shape space.

\end{itemize}
Note that the third point leads to more efficient Algorithms in cases where one only needs a shape geodesic and not the optimal registration between surfaces. 
It provides the same geodesic path despite arbitrary initial parameterizations (or registrations) of given surfaces, 
and saves the computational cost of finding a registration.
This fact is also a source of limitation 
in the situation where one needs a registration. If one wants to use geodesic lengths for comparing shapes, then a registration is not needed. 
However, if one wants to study statistical summaries of deformation fields, then a registration will be needed.

The rest of this paper is organized as follows. 
Section 2 describes the mathematical representation of embedded surfaces and establishes mathematical setup. Section 3 is devoted to the description of gauge invariance and to the definition of the Riemannian metric involved in this paper. The geodesic computation is 
described in Section 4 and Section 5 presents the experimental results.

\section{Mathematical Setup}

\subsection{Notation}
We will represent a shape $S$ with an embedding $f~:\mathbb{S}^2\rightarrow\mathbb{R}^3$ 
such that the image $f(\mathbb{S}^2)$ is $S$. The function $f$ is also called a \textit{parameterization} of the surface $S$. 

We will use local coordinates $(u, v)$ on the sphere. For the theoretical framework, any coordinates on the sphere are suitable, but in the application we use spherical coordinates~: $u$ stands for the polar angle and ranges from $0$ to $\pi$, and $v$ denotes azimuthal angle and ranges from $0$ to $2\pi$.

Recall that a map $f~:\mathbb{S}^2\rightarrow\mathbb{R}^3$ is an \textit{embedding} when: for any point $(u,v) \in \mathbb{S}^2$,
(1)  $f$ is smooth, in particular the derivatives $f_u$ and $f_v$ of $f$ with respect to $u$ and $v$ are well-defined,
(2)  $f$ is an immersion, i.e. the cross product $f_u \times f_v$ never vanishes and allows us to define the normal  (resp. tangent) space to the surface  
$f(\mathbb{S}^2)$ at a point $f(u, v)$ as the subspace of $\mathbb{R}^3$ which is generated by (resp. orthogonal
to) $f_u \times f_v$, and (3) $f$ is an homeomorphism onto its image, i.e. points on $f(\mathbb{S}^2)$ that look close
in $\mathbb{R}^3$ are images of close points in $\mathbb{S}^2$.
If $f$ is an embedding,  
then the surface $f(\mathbb{S}^2)$ is naturally oriented by the frame $\{f_u, f_v\}$, or equivalently by the normal vector field $f_u \times f_v$.

We define the space of all such surfaces as
$$\mathcal{F} := \{f~:\mathbb{S}^2\rightarrow \mathbb{R}^3, f \textrm{~is~an~embedding}\}.$$
It is often called the \textit{pre-shape space} since objects with same shape but 
different orientations or parameterizations may correspond to different points in $\mathcal{F}$. 
The set $\mathcal{F}$ is itself a manifold, as an open subset of the linear space $\mathcal{C}^{\infty}(\mathbb{S}^2, \mathbb{R}^3)$ of smooth functions from 
$\mathbb{S}^2$ to $\mathbb{R}^3$ (see Theorem~3.1 in \cite{BBM} and the references therein). The tangent space to $\mathcal{F}$ at $f$, denoted by $T_f\mathcal{F}$, is therefore just $\mathcal{C}^{\infty}(\mathbb{S}^2, \mathbb{R}^3)$.

The shape-preserving transformations of 3D object can be expressed as group actions on $\mathcal{F}$.
The group $\mathbb{R}^+$ with multiplication operation acts on $\mathcal{F}$
by \textit{scaling}~: $\left(\beta, f\right) \mapsto \beta f$, for $\beta\in\mathbb{R}^+$ and $f\in \mathcal{F}$.
The group $\mathbb{R}^3$ with addition as group operation acts on $\mathcal{F}$, 
by \textit{translations}~: $(v, f) \mapsto f+v$, for $v\in \mathbb{R}^3$ and $f\in \mathcal{F}$.
The group $\textrm{SO}(3)$ with matrix multiplication as group operation  acts on $\mathcal{F}$, 
 by \textit{rotations}~: $(O, f)\mapsto O f$, for $O\in \textrm{SO}(3)$ and $f\in \mathcal{F}$. 
Finally, 
the group $\Gamma := \textrm{Diff}^+(\mathbb{S}^2)$ consisting of diffeomorphisms which preserve the orientation of $\mathbb{S}^2$ acts also on $\mathcal{F}$,
by \textit{reparameterization}~: $(\gamma, f)\mapsto f\circ\gamma^{-1}$, for $\gamma\in \textrm{Diff}^+(\mathbb{S}^2)$ and $f\in\mathcal{F}$.
The use of $\gamma^{-1}$, instead of $\gamma$, ensures that the action is from left and, since the action of $\textrm{SO}(3)$ is also 
from left, one can form a joint action of $G := \textrm{Diff}^+(\mathbb{S}^2)\times  \textrm{SO}(3)\rtimes \mathbb{R}^3$ on $\mathcal{F}$.
In this paper, the translation group is taken care of by using a translation-independant metric (the elastic metric) and, when needed, 
the scaling is taken care of by rescaling the surfaces to have unit surface area. Therefore, in the following we will 
focus only on the reparameterization group $\Gamma$ and on the rotation group $\textrm{SO}(3)$.

\subsection{Shape Space as quotient space}
Since we are only interested in shapes of surfaces, 
we would like to identify surfaces that can be related through a shape-preserving 
transformation. This is accomplished using the notion of group action and orbits
under those group actions.

Given a group $G$ acting on $\mathcal{F}$, the elements in $\mathcal{F}$ obtained by following a fix 
parameterized surface $f\in\mathcal{F}$ when 
acted on by all elements of $G$ is called the $G$-\textit{orbit} of $f$ or the \textit{equivalence class} of $f$ under the action of $G$,
and will be denoted by $[f]$. 
In particular, when $G$ is the reparameterization group, the orbit of $f\in \mathcal{F}$ is characterized by the surface $f(\mathbb{S}^2) = S$, 
i.e. the elements in $[f] = \{f\circ \gamma^{-1} \textrm{~for~}\gamma \in \Gamma\}$ are all possible parameterizations of $S$. For instance in Fig. \ref{fig_illustration}, the first column contains some parameterized horses that are elements of the same orbit.  
The set of orbits of $\mathcal{F}$ under a group $G$ is called the \textit{quotient space} and will be denoted by $\mathcal{F}/G$. 
The quotient space of interest in this paper is called \textit{shape space} and is defined as follows.
\begin{definition}
The \textit{shape space}  $\mathcal{S}$ is the set of oriented surfaces in $\mathbb{R}^3$, which are diffeomorphic to $\mathbb{S}^2$, 
modulo translation and rotation. It is isomorphic to the quotient space of the pre-shape space  
$\mathcal{F}$ by the shape-preserving group $G := \textrm{Diff}^+(\mathbb{S}^2)\times \textrm{SO}(3)\rtimes \mathbb{R}^3$~:
$\mathcal{S} = \mathcal{F}/G$. 
\end{definition}

It is important to note that
the shape space $\mathcal{S} = \mathcal{F}/G$ is a smooth manifold and the 
canonical projection $\Pi~:\mathcal{F}\rightarrow \mathcal{F}/G$, $f\mapsto [f]$ is a submersion (see for instance \cite{BE} and \cite{M}).
This submersion is useful in establishing the notion of a vertical space that will be needed a little later. 
By definition, the vertical space of a submersion  is the kernel space of its differential. When the submersion is a quotient map by a group action, the vertical space is the tangent space to the orbit (the terminology comes from the fact that the orbits are usually depicted as vertical fibers over a base manifold which is the quotient space, see Fig.\ref{fig_illustration}).
In the case of the submersion $\tilde{\Pi}: \mathcal{F}\mapsto \mathcal{F}/\textrm{Diff}^+(\mathbb{S}^2)$, the vertical space
takes a very natural, intuitive form. 
\begin{proposition}
The vertical space $Ver(f)$ of $\tilde{\Pi}$ at some embedding $f\in\mathcal{F}$ is the space of vector fields which are tangent to the 
shape $f(\mathbb{S}^2)$, or equivalently the space of vector fields such that the dot product with the unit normal vector field $n_f := \frac{f_u \times f_v}{\|f_u \times f_v\|}~:\mathbb{S}^2 \rightarrow\mathbb{R}^3$ vanishes~:
$$  
Ver(f) =  \{\delta f:\mathbb{S}^2\rightarrow \mathbb{R}^3| \delta f(s)\cdot n_f(s) = 0, \forall s \in \mathbb{S}^2\}.
$$
\end{proposition}
\begin{remark}
A canonical complement to this vertical space (consisting of tangent vector fields) 
is given by the space of vector fields normal to the surface $f(\mathbb{S}^2)$
denoted by $Nor$. This is the sub-bundle of the tangent bundle $T\mathcal{F}$ defined by
$$
Nor(f) = \{\delta f:\mathbb{S}^2\rightarrow \mathbb{R}^3| \delta f(s)\times n_f(s) = 0,\forall s \in\mathbb{S}^2\}.
$$
Any tangent vector $\delta f\in T_f\mathcal{F}$ admits a unique decomposition 
$\delta f = \delta f^{T} + \delta f^{\perp}$
into its tangential part $\delta f^{T}\in Ver(f)$ and its normal part $\delta f^{\perp}\in Nor(f)$. 
Specifically, the normal part is given by: 
\begin{equation}\label{eq:normal-component}
\delta f^{\perp} =  \left(\delta f\cdot n_f\right) n_f
\ .
\end{equation}
See Fig. \ref{vector_decomposition} for an illustration of this decomposition.
Generally speaking,  one has $T\mathcal{F} = Ver \oplus Nor$
as a direct sum of smooth fiber bundles over $\mathcal{F}$.
This decomposition is preserved by the action of the reparameterization group $\Gamma$, i.e. $\left(\delta f\circ\gamma\right)^{T} = 
 \delta f^{T} \circ\gamma$ and $\left( \delta f\circ\gamma\right)^{\perp} =  \delta f^{\perp}\circ\gamma$  (for a proof of this statement, see Section 1 of the Supplementary Material).\color{black} 
 \end{remark}
 
 The interest in splitting a perturbation $\delta f$ into its normal and vertical components comes 
 from the fact that the vertical component $\delta f^T \in Ver(f)$ can only lead to a shape-preserving transformations of the surface $f(\mathbb{S}^2)$. 
 Thus, in the process of deforming one shape into another (for instance along a geodesic path) and
 quantifying shape differences between them using geodesic lengths, we are not 
 interested in measuring deformations that are in $Ver(f)$. An important novelty of this paper is that the
eventual Riemannian metric is imposed only on the $\delta f^{\perp}$ components of the 
perturbations, and that the  $\delta f^T$ components have a zero contribution to the metric.

 \section{Gauge Invariance and Riemannian Metric}
 As mentioned earlier, another important goal of this paper is in developing a framework 
 that is gauge invariant. To appreciate the utility of this framework, we first provide a 
 precise definition and then motivate its use in shape analysis. 
 
 \subsection{Defining Gauge Invariance}
 The gauge invariance relates to the parameterization of surfaces along a path in $\mathcal{F}$ and, thus, 
the mathematical objects of importance in this section 
are paths $\Psi~:[0, 1]\mapsto \mathcal{F}$. The set of such paths is the smooth manifold 
$\mathcal{P} :=\mathcal{C}^{\infty}([0,1], \mathcal{F})$. 

An element of $\mathcal{P}$ can be thought of as a metamorphosis from the initial shape to the final shape. For instance,
Fig. \ref{fig_illustration} shows two elements in $\mathcal{P}$
as two different deformations from a parameterized horse to a parameterized cat. To have a picture in mind, consider the upper 
path $\Psi~: t\mapsto \Psi(t)$ in $\mathcal{P}$: at each time step $t\in  [0, 1]$, $\Psi(t)$ is a parameterized shape, i.e. a map from our model manifold $\mathbb{S}^2$ into $\mathbb{R}^3$. 
The map $\Psi(0)$ is the parameterization of our initial parameterized shape chosen to be a horse and $\Psi(1)$ is the parameterization of our final parameterized shape which, in this case, is a cat.

The definition of length of the path $t\mapsto \Psi(t)$ requires specification of a metric on $\mathcal{F}$. Given such 
a metric $(\!(\cdot, \cdot )\!)$,
one can define the length as:
\begin{equation}\label{length}
L[\Psi] = \int_{0}^{1} (\!( \Psi_t(t), \Psi_t(t))\!)_{\Psi(t)}^{\frac{1}{2}} dt,
\end{equation}
where $\Psi_t(t)=\frac{d\Psi}{dt}(t)$ is the velocity vector of the path $t\mapsto \Psi(t)$, i.e. an infinitesimal deformation of the parameterized shape $\Psi(t)$.
The geodesic distance between two shapes $f_1$ and $f_2$ is then defined by
\begin{equation}\label{distance}
d(f_1, f_2) = \inf_{\Psi:[0,1] \to \mathcal{F}| {\Psi(0) = f_1, \Psi(1) = f_2}}~L[\Psi],
\end{equation}
where the infimum is taken over all paths connecting shape $f_1$ and shape $f_2$.

We would like the length $L[\Psi]$, for any path $\Psi$,  
 to match the length of the path $t\mapsto \Psi(t)\circ \gamma(t)$, where $t\mapsto \gamma(t) \in \Gamma$ is any time-dependent reparameterization of $\mathbb{S}^2$~:
\begin{equation}\label{gauge_invariance}
L[\Psi] = L[\tilde{\Psi}],\ \ \mbox{where}\ \tilde{\Psi}(t) = \Psi(t)\circ\gamma(t).
\end{equation}
More formally, set $\Gamma = \textrm{Diff}^+(\mathbb{S}^2)$ and define the group 
$   \mathcal{G}~:= \mathcal{C}^{\infty}([0,1], \Gamma)$,
    of time-dependant reparameterizations that acts on $\mathcal{P}$  according to
   $$
   \begin{array}{clc}
   \mathcal{G}\times\mathcal{P}&\longrightarrow&\mathcal{P}\\
   (t\mapsto \gamma(t), t\mapsto \Psi(t))&\longmapsto &(t\mapsto \Psi(t)\circ\gamma(t)).
   \end{array}
   $$
The group $\mathcal{G}$  is called the \textit{gauge group}, and one says that $\mathcal{G}$ acts by \textit{gauge transformations}. 
We are looking for a framework where the length of a path is \textit{invariant to gauge transformations}, i.e. satisfies Eqn.~\eqref{gauge_invariance}.   
One should distinguish these transformations from {\it temporal} reparameterizations of the path $\Psi$ itself. A gauge transformation changes 
spatial reparameterization of surfaces, while preserving shapes, along the path, while a temporal reparameterization changes 
the time it takes to reach each shape along the path.

To build a gauge invariant framework, the basic idea is as follows:  
take any $\Gamma$-invariant Riemannian metric  
$\langle\!\langle \cdot, \cdot \rangle\!\rangle$ on the pre-shape space, 
and \textit{ignore} the direction tangent to the reparameterization orbit. 
(An example of $\Gamma$-invariant Riemannian metric is the elastic metric defined in Eqn.~\eqref{eqn:elastic-metric} 
as is shown in Section 2 of the Supplementary Material). 
More precisely,
let  $\langle\!\langle\cdot, \cdot\rangle\!\rangle$ be a Riemannian metric on pre-shape space $\mathcal{F}$ which is  preserved by the action of the group of reparameterizations $\Gamma$, that is:
 \begin{equation}\label{condition1}
 \langle\!\langle  \delta f_1\circ\gamma,  \delta f_2\circ\gamma \rangle\!\rangle_{f\circ\gamma} = \langle\!\langle  \delta f_1,  \delta f_2\rangle\!\rangle_{f},
 \end{equation}
for any  $f \in\mathcal{F}$, for any $ \delta f_1, \delta f_2 \in T_{f}\mathcal{F}$ and any $ \gamma\in\Gamma$. Given a $\Gamma$-invariant sub-bundle $H$ of $T\mathcal{F}$ 
such that 
\begin{equation}\label{direct_sum}
H(f) \oplus Ver(f) = T_{f}\mathcal{F},
\end{equation} 
denote by  $p_{H}: T_{f}\mathcal{F}\rightarrow H(f)$ the projection onto $H(f)$ with respect to the direct sum decomposition given in Eqn.~\eqref{direct_sum}. 
 \begin{figure}[ht]
 		\centering
 		\includegraphics[width=7cm]{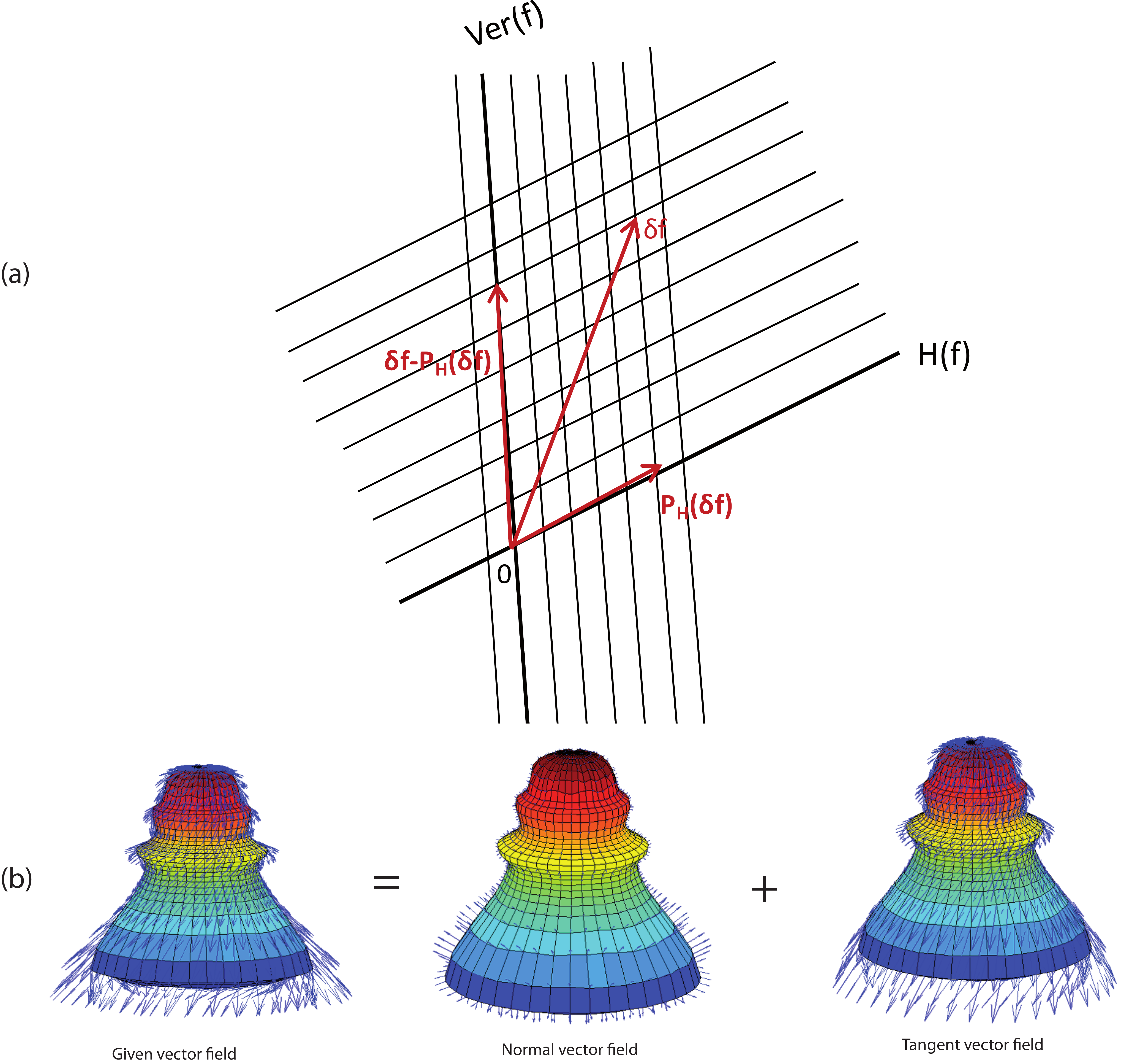}
 		\caption{\scriptsize a. Direct sum decomposition $H(f) \oplus Ver(f) = T_{f}\mathcal{F}$. 
        b.~Vector field decomposition into tangent and normal directions}
		\label{vector_decomposition}
 		\end{figure}
This means that any element $\delta f\in T_f\mathcal{F}$ admits a unique decomposition into the sum of an element  $p_H(\delta f)$ in $H(f)$ 
and an element in $Ver(f)$. 
We illustrate this decomposition of vector spaces in Fig.~\ref{vector_decomposition}.a, while
the particular case when $H$ is the space of normal vector fields $Nor$ is shown in Fig.~\ref{vector_decomposition}.b.
\begin{proposition}\label{theorem}
The non-negative semi-definite inner product on pre-shape space defined by
$$
\left( \!\left( \delta f_1,  \delta f_2\right)\!\right)_f:= \langle\!\langle  p_H(\delta f_1),  p_H(\delta f_2)\rangle\!\rangle_{f}
$$
satisfies the gauge-invariance condition given in Eqn.~\eqref{gauge_invariance} and 
 induces a Riemannian metric on quotient space $\mathcal{S}$ such that the quotient map is an isometry between $H(f)$ and the tangent space $T_{[f]}\mathcal{S}$.
 \end{proposition}
 
\subsection{Distinction between Gauge Invariant Framework and Quotient Riemannian Framework}
In practice the subbundle $H$ has to be chosen in order to make the implementation easy.
A natural choice of subbundle $H$ is the normal bundle $Nor$ which is preserved by the action of the reparameterization group $\Gamma$  
(for a proof of this statement, see Section 1 of the Supplementary Material). 
We have used this subbundle in the present paper. Another requirement is that the chosen Riemannian metric has to be $\Gamma$-invariant. 
This is the case for the elastic metric defined in next section. 
We will therefore apply the idea of gauge invariance to the concrete example of the elastic metric and the normal bundle $Nor$ in the remainder of this paper.
It is worth noting that the Riemannian metric on shape space obtained by restricting a Riemannian metric on preshape space to the normal bundle $Nor$ differs in general from the quotient Riemannian metric. In fact, the quotient metric coincides with the restriction to the subbundle $Nor$ if and only if the \textit{Horizontal} subbundle defined by     $Hor(f) = \textrm{Ker}(d\pi)^{\perp}$ is the normal bundle. This is not the case for the elastic metric. 
We also remark that the present gauge invariant framework has been used implicitly in \cite{BHM3}, Section~6, and \cite{BHM5}, Section~11, in the case where the horizontal bundle coincides with the normal bundle.

 \subsection{Elastic Riemannian Metric}
 Next, we will choose a Riemmanian metric on $\mathcal{F}$ that will enable a gauge-invariant analysis as
 stated above. We will use the elastic Riemannian metric proposed by Jermyn et al \cite{JermynECCV2012} and given in 
 Eqns.~\eqref{eqn:elastic-metric} and \eqref{eqn:elastic-metric2}. However, before we use this metric we motivate its use 
 by making a connection between the space of parameterized surfaces $\mathcal{F}$ and the
 space of metrics on a domain, and we will provide some geometrical interpretation 
 of terms in that elastic  metric. 
The space of positive-definite Riemannian metrics on 
$\mathbb{S}^2$ will be denoted by $\operatorname{Met}(\mathbb{S}^2)$. 
Consider a parameterized surface $f~: \mathbb{S}^2 \rightarrow \mathbb{R}^3$. Denote by $g = f^*\bar{g}$ the pull-back of the Euclidian metric $\bar{g}$ of 
$\mathbb{R}^3$
and by $n_f$ the unit normal vector field (Gauss map) on $S = f(\mathbb{S}^2)$. 

The metric $g$ and the normal vector field $n_f$ are defined using derivatives of $f$ according to:
\begin{eqnarray*}
g &=&\left( \begin{smallmatrix} f_u\cdot f_u & 
f_u\cdot f_v \\ 
f_v\cdot f_u & f_v\cdot f_v
\end{smallmatrix} \right)= \operatorname{Jac}(f)^{T}\operatorname{Jac}(f),\\
&=& \left( \begin{smallmatrix} E & F\\ 
F & G \end{smallmatrix}  \right),\ \  \operatorname{Jac}(f) = [f_u  \,\,f_v], \textrm{~and}\\
n_f &=& \frac{f_u\times f_v}{\|  f_u\times f_v\|}
,\ \ \left\| f_u\times f_v\right\| = \sqrt{\det{g}}\ = |g|^{\frac{1}{2}},
\end{eqnarray*}
where $f_u$ and $f_v$ are the derivatives of $f$ with respect to the local coordinates  $(u,v)$  on the sphere.
We consider the following relationship between parameterized
surfaces on one hand and the product space of metrics and normals on the other~:
$$
\begin{array}{lclc}
\Phi~:& \mathcal{F} &\longrightarrow &\operatorname{Met}(\mathbb{S}^2) \times \mathcal{C}^{\infty}(\mathbb{S}^2, \mathbb{S}^2)\\
& f & \longmapsto &  (g,  {n_f}).
\end{array}
$$

It follows from  the fundamental theorem of surface theory (see Bonnet's Theorem in \cite{DoCarmo} for the local result, Theorem 3.8.8 in  \cite{Klingenberg} or Theorem 2.8-1 in \cite{Ciarlet}  
for the global result) that two parameterized surfaces $f_1$ and $f_2$ having the same representation $(g, n)$ differ at most by a translation and 
rotation. This is an important result, and implies that we can represent a surface 
by its induced metric $g = f^*\bar{g}$ and the unit normal field $n = n_f$, for the purpose of analyzing its shape. 
We will not loose any information about the shape of a surface $f$ if we represent
it by the pair $(g,n)$. 
Let $\delta f_1$, $\delta f_2$ denote two perturbations of a surface $f$, and 
let $(\delta\! g_1, \delta\!  {n}_1) = \Phi_*(\delta f_1)$, $(\delta\! g_2, \delta\!  {n}_2) = \Phi_*(\delta f_2)$ denote
the corresponding perturbations in $(g,n)$ of $f$. The expression for $\Phi_*$ is given by: 
\begin{eqnarray*}
\delta\! g  &=& \operatorname{Jac}(f)^T \operatorname{Jac}(\delta\!f) + (\operatorname{Jac}(\delta\!f))^T \operatorname{Jac}(f) \\
&=&   \left(\begin{array}{cc}2 f_u\cdot \delta\!f_u & 
f_u\cdot \delta\!f_v + f_v\cdot \delta\!f_u  \\ 
f_u\cdot \delta\!f_v +  f_v\cdot \delta\!f_u  
&  2 f_v\cdot \delta\!f_v \end{array}  \right),
\\
\delta\! n &=& -\frac{1}{2}\operatorname{Tr}(g^{-1} \delta\! g)\, n+ \frac{1}{|g|^{\frac{1}{2}}}
\left( \delta\!f_u\times  f_v +  f_u\times   \delta\!f_v  \right).
 \end{eqnarray*}

Then, by definition, the metric on $\mathcal{F}$ used in the present paper measures these perturbations using the expression
\begin{eqnarray}\label{eqn:elastic-metric}
\langle\!\langle \delta f_1, \delta f_2 \rangle\!\rangle_{f} &=& 
 \int_{\mathbb{S}^2} ds |g|^{\frac{1}{2}} \left\{ a \operatorname{Tr}(g^{-1}\delta\! g_1 \,g^{-1} \delta\! g_2) \right. \nonumber \\
&& \hspace*{-0.5in}  +
\frac{\lambda}{2} \operatorname{Tr}(g^{-1}\delta\! g_1)\operatorname{Tr}(g^{-1}\delta\! g_{2}) \left.+ c \delta\!  {n}_1\cdot \delta\!  {n}_2 \right\}.
\end{eqnarray}
The same metric (with $a=1$) was introduced in \cite{JermynECCV2012}, Eqn.~(2), and called ``elastic metric''. 
A related metric measuring the elastic deformation of the interiors of shapes was used in \cite{FJSY09} (see Eqn.~(4) in \cite{FJSY09}).
The metric given in Eqn.~\eqref{eqn:elastic-metric} can be decomposed into three parts
 \begin{eqnarray}\label{eqn:elastic-metric2}
 \langle\!\langle \delta f_1, \delta f_2  \rangle\!\rangle_{f}  
 &=&  
 \hspace*{-0.1in} \int_{\mathbb{S}^2} ds |g|^{\frac{1}{2}} \left\{ a\operatorname{Tr} \left( (g^{-1}\delta\! g_1)_{\textbf{0}} \,(g^{-1} \delta\! g_2)_{\textbf{0}}\right) \right.
 \nonumber \\ 
  && \hspace*{-0.5in} \left. + 
b \operatorname{Tr}(g^{-1}\delta\! g_1)\operatorname{Tr}(g^{-1}\delta\! g_{2})  + c \delta\!  {n}_1\cdot \delta\!  {n}_2 \right\},
 \end{eqnarray}
 where $b = \frac{\lambda+a}{2}$ and where $A_{\textbf{0}} $ is the traceless part of a $2\times 2$-matrix $A$ defined as $A_{\textbf{0}}  = A -\frac{\operatorname{Tr}(A)}{2}I_{2\times2}$. 
 The term multiplied by $a$ measures area-preserving changes in the induced metric $g$, the term multiplied by $b$ measures changes in the area of patches, and the last term measures bending. 
 Note that only the relative weights $b/a$ and $c/a$ are meaningful.

Now we consider a key property of this metric that relates to reparameterization of surfaces. 
Recall that $\Gamma := \operatorname{Diff}^+(\mathbb{S}^2)$ denotes the subgroup of $\operatorname{Diff}(\mathbb{S}^2)$ consisting of diffeomorphisms $\gamma$ which preserve the orientation of $\mathbb{S}^2$, i.e. such that
$\det \operatorname{Jac}(\gamma) >0$. 
(Note that for a diffeomorphism $\gamma\in \operatorname{Diff}(\mathbb{S}^2)$, 
since $\operatorname{Jac}(\gamma)$ is invertible, the determinant of $\operatorname{Jac}(\gamma)$ never vanish. It follows that either 
$\det \operatorname{Jac}(\gamma)(s)>0$ for all $s\in \mathbb{S}^2$, or $\det \operatorname{Jac}(\gamma)(s)<0$ for all $s\in \mathbb{S}^2$.)
It will be called the group of orientation-preserving reparameterizations.
The group 
$\Gamma = \operatorname{Diff}^+(\mathbb{S}^2)$  acts on $\operatorname{Maps}(\mathbb{S}^2, \mathbb{R}^3)$ by pre-composition.
That is, a surface $f$ is reparameterized by a $\gamma \in \operatorname{Diff}^+(\mathbb{S}^2)$ according to 
$f \mapsto f \circ \gamma^{-1}$. How does the metric-normal representation $(g,n)$ of that surface change due to 
reparameterization? This representation of the reparameterized surface is given by $(\gamma^{-1*}g,  {n}\circ \gamma^{-1})$.
This representation is $\Gamma$-equivariant for the actions introduced, i.e. if we reparameterize a surface and then compute
its $(g,n)$ representation, or if we compute $(g,n)$ representation of a surface and then reparameterize 
them according to $(\gamma^*g,  {n}\circ \gamma)$, we get the same result.

\begin{proposition}\label{Gamma-invariance}
The elastic metric is  invariant to the action of $\operatorname{Diff}^+(\mathbb{S}^2)$. 
\end{proposition}
{\bf Proof}: Please refer to  Section 2 of the Supplementary Material. 

Although this elastic metric has been introduced by Jermyn et al. \cite{JermynECCV2012}, it has not been used completely
for shape analysis of surfaces. Furthermore, we are going to use it in a novel way -- by 
restricting its evaluation only to the normal vector fields on a surface (see next section for a geometric expression of the resulting metric on shape space). 

\begin{definition}
\label{def_metric}
For any two perturbations $\delta f_1, \delta f_2 \in T_f\mathcal{F}$ define the pairing
$$
\left(\!\left(\delta f_1,\delta f_2\right)\!\right)_{f} = 
\innerd{\delta f_1^{\perp} }{\delta f_2^{\perp} }_f\ ,
$$
where $\delta f_i^{\perp}$ is the normal component of $\delta f_i$ as defined in Eqn.~\eqref{eq:normal-component} and where  $\innerd{\cdot}{\cdot}_f$ is as given in 
Eqn.~\eqref{eqn:elastic-metric2}. 
\end{definition}

\begin{remark}
It follows from proposition \ref{theorem}, that $\left(\!\left(\cdot,\cdot\right)\!\right)$ satisfies the gauge-invariant condition 
$L[\Psi] = L[\tilde{\Psi}]$,
where $\Psi$ is any path of shapes, $\tilde{\Psi}(t) = \Psi(t)\circ\gamma(t) $ with $t\mapsto \gamma(t)$ any time-dependant reparameterization, and
$L[\Psi]$ is as specified in Eqn.~\eqref{length}. 
\end{remark}

 \subsection{Geometric expression of the elastic metric in the normal direction} \label{sec:geometric-interp}
In this section,  we will give some geometric interpretation of the restriction of the elastic metric on the space of normal vector fields introduced in 
the previous section. Given a surface $f$ parameterized by $(u,v)$, 
we will consider normal variations:
$f_{\varepsilon}(u, v) = f(u, v) + \varepsilon h(u, v) n(u, v)$,
where $(u, v)\in \mathbb{S}^2$,  $\varepsilon > 0$, $n(u,v) = n_f(u,v)$ is the unit normal to the surface $f(\mathbb{S}^2)$ at 
$f(u, v)$, and $h~: \mathbb{S}^2 \rightarrow\mathbb{R}$ is a real function corresponding to  the amplitude of the normal vector field $h\,n_f$.
Let us compute the first fundamental form $g_{\varepsilon}$ of the surface parameterized by  $f_{\varepsilon}$, i.e. 
the metric induced on  the parameterized surface $f_{\varepsilon}$ by the Euclidian metric of $\mathbb{R}^3$. We obtain
\begin{equation}\label{f_epsilon}
\begin{array}{l}
f_{\varepsilon, u}:= \frac{\partial f_{\varepsilon}}{\partial u} = f_{u} + \varepsilon h n_{u} + \varepsilon h_{u} n,\\
\\
f_{\varepsilon, v}:= \frac{\partial f_{\varepsilon}}{\partial v} = f_{v} + \varepsilon h n_{v} + \varepsilon h_{v} n.
\end{array}
\end{equation}
Therefore
$$
f_{\varepsilon, u}\cdot f_{\varepsilon, u} = f_u\cdot f_u  + 2 \varepsilon h  n_u\cdot f_u + \varepsilon^2\left(h^2  n_u\cdot n_u  + h_u^2\right),
$$
where we have used that $n\cdot f_u  =  0 $ and $n_u\cdot n  = 0$ since $n\cdot n  = 1$.
Similarly
$$
 f_{\varepsilon, v}\cdot f_{\varepsilon, v} = f_v\cdot  f_v  + 2 \varepsilon h n_v \cdot f_v + \varepsilon^2\left(h^2  n_v\cdot n_v  + h_v^2\right),
$$
and 
$$
 f_{\varepsilon, u}\cdot f_{\varepsilon, v} =  f_u\cdot f_v  +  \varepsilon h \left( n_u\!\cdot\! f_v \!+\!  f_u\!\cdot\! n_v\right) +
\varepsilon^2\left(h^2  n_u\!\cdot\! n_v \!+\! h_u h_v\right).
$$
It follows that
$$
\begin{array}{ll}
g_{\varepsilon} = &  g 
+ 2\varepsilon h \left(\begin{array}{cc}n_u\cdot f_u &  n_u\cdot f_v\\  n_v\cdot f_u &  n_v\cdot f_v\end{array}\right) \\ & 
\!\!\!\!\!\!+ \varepsilon^2 h^2 \left(\begin{array}{cc}  n_u\cdot n_u & n_u\cdot n_v\\  n_v\cdot n_u &  n_v\cdot n_v\end{array}\right)
 + \varepsilon^2 \left(\begin{array}{cc} h_u^2 & h_u h_v\\ h_u h_v & h_v^2 \end{array}\right).
 \end{array}
$$
Using the definition of the second fundamental form $\textbf{II}$ of the surface $f(\mathbb{S}^2)$, we obtain
$$
g_{\varepsilon} = g -2\varepsilon h \textbf{II} + \varepsilon^2 h^2 \textbf{II}g^{-1} \textbf{II}+ \varepsilon^2 \left(\begin{array}{cc} h_u  & h_v\end{array}\right)^T\left(\begin{array}{cc} h_u  & h_v\end{array}\right)\!.
$$
It follows that 
\begin{equation}\label{g_shape_operator}
g^{-1} \delta g = -2 h g^{-1} \textbf{II} = -2 h \textrm{\textbf{L}},
\end{equation}
where $\textrm{\textbf{L}}$ is called the {\it shape operator}. 
Recall that the eigenvalues of $\textrm{\textbf{L}}$ are the principal curvatures of the surface $f(\mathbb{S}^2)$, denoted by $\kappa_1$ and $\kappa_2$, which provide local information about the surface: at a given point on the surface, they measure the greatest and smallest possible curvatures of a curve drawn on the surface passing through this point. For instance, the vanishing of the principal curvatures at one point of the surface tells that the surface is flat near this point (i.e. looks like a plane). The equality  $\kappa_1=\kappa_2=1/R$ at one point tells that the surface looks like a sphere of radius $R$ near this point. In other words, $\kappa_1$ and $\kappa_2$ are functions on the surface that characterize how the surface is locally curved.

%

On the other hand, the variation $\delta\!n$ of the normal vector field satisfies 
$\delta\!n\cdot n = 0$ since the norm of $n$ remains constant. Moreover $n \cdot f_u =  n\cdot f_v = 0$, 
therefore $\delta\!n\cdot f_u = -  n\cdot \delta\!f_u$ and
$\delta\!n\cdot f_v = - n\cdot \delta\!f_v$.
By Eqn.~\eqref{f_epsilon}, $\delta\!f_u = h n_u + h_u n$, hence $\delta n\cdot f_u = -h_u$ and similarly $ \delta n\cdot f_v = -h_v$.
Consequently $\delta\!n = \alpha f_u + \beta f_v$ where 
$
\left(\begin{smallmatrix} \alpha\\ \beta\end{smallmatrix}\right) = -g^{-1} \left(\begin{smallmatrix} h_u\\ h_v\end{smallmatrix}\right).
$
It follows that for two normal vector fields
$hn$ and $kn$ with $h, k\in\mathcal{C}^{\infty}(\mathbb{S}^2, \mathbb{R})$, one has
\begin{equation}\label{delta_epsilon_n}
 \delta\!n_1\cdot \delta\!n_2  =
\left(\begin{smallmatrix} h_u & h_v\end{smallmatrix}\right) g^{-1} \left(\begin{smallmatrix} k_u\\k_v\end{smallmatrix}\right).
\end{equation}
Using 
Eqn.~\eqref{g_shape_operator}  and Eqn.~\eqref{delta_epsilon_n} the 
elastic metric restricted to these normal fields is given by~:
\begin{eqnarray}\label{expression_metric}
 \left(\!\left( hn, kn\right)\!\right)_{f}&\!\!\!\!\!\!\!\!\!\!\!\!\!\!\!\!=  \int_{\mathbb{S}^2} ds |g|^{\frac{1}{2}} \left\{hk \left(2a (\kappa_1 - \kappa_2)^2   \right.  \right.\nonumber\\
&\left. + 4b (\kappa_1+\kappa_2)^2\right) \!+\! \left. c  \left(\begin{smallmatrix} h_u & h_v\end{smallmatrix}\right) g^{-1}\!\left(\begin{smallmatrix} k_u\\k_v\end{smallmatrix}\right)\right\}\!\!.
\end{eqnarray}
This is the form used to define and compute geodesic paths in the shape space $\mathcal{S}$ in this paper. 
The difference $\kappa_1 - \kappa_2$ in the first  term has been called the normal deformation of the surface in \cite{Ivanova}. 
The sum $\kappa_1 + \kappa_2$ is twice the mean curvature which measures variations of the area of local patches. These two terms are related to the shape index $\textrm{idx} = \frac{2}{\pi} \textrm{arctan}\frac{\kappa_1+ \kappa_2}{\kappa_1-\kappa_2}$ \cite{koenderink1990solid}. The last term in Eqn.~\eqref{expression_metric} measures variations of the normal vector field, i.e. bending.
%

%
%
%
%

\section{Geodesic Computation}


Finding geodesics between two surfaces $f_1$ and $f_2$ under invariant Riemannian metrics is a difficult problem. 
In the present case, analytical solutions are not known and
we will use a path-straightening approach to find geodesics. This method has been used for instance 
in \cite{KurtekPAMI2012} and \cite{JermynECCV2012}. The basic idea here is to connect $f_1$ and $f_2$ by any 
initial path and then iteratively straighten it until it becomes a geodesic. The update is performed using the gradient of an energy function. 
As mentioned earlier, this method only achieves a local minimum of the energy function, resulting in a geodesic path that 
may not be the shortest geodesic.

\subsection{Removing rotations and translations}\label{rotation}

Since we are only interested in shapes of objects and not in the way objects are oriented or placed  in the ambient space $\mathbb{R}^3$, 
we have to remove the actions of the rotation and translation groups.
In theory, we could deal with them as we do with the group of reparameterizations. 
However,   since $\textrm{SO}(3)\rtimes \mathbb{R}^3$ is just a $6$-dimensional Lie group in comparison to the 
infinite-dimensional Fr\'echet Lie group $\textrm{Diff}^+(\mathbb{S}^2)$, it is more efficient to 
do the following. First find the best translation and rotation 
that align the two objects to be compared, and then find the geodesic between them. 
To center an object we use Algorithm~3 given in Section~3 of the Supplementary Material to compute the center of mass and then subtract it from the surface coordinates.

\begin{figure}[!ht]
 		\centering
 		\includegraphics[width=8cm]{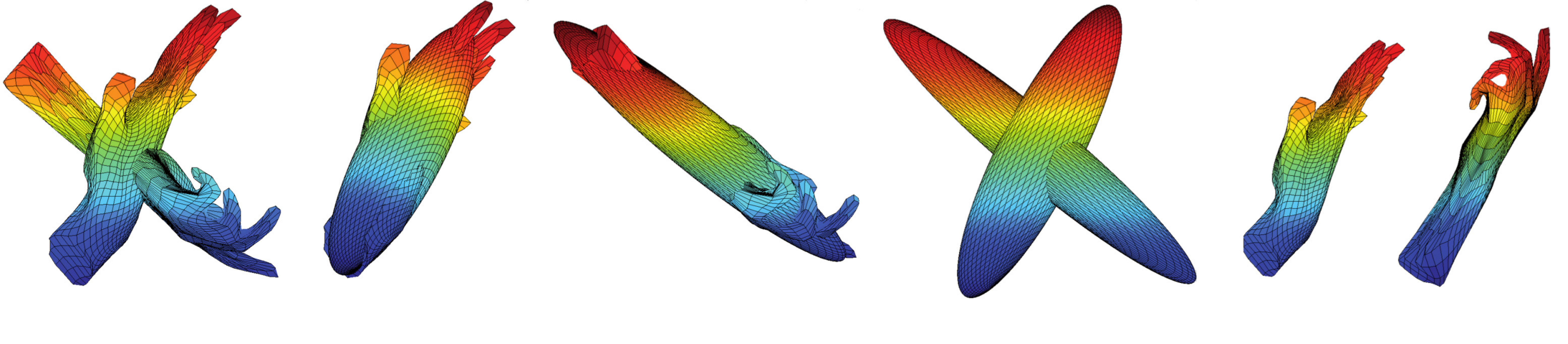}
 		\caption{\scriptsize Rotational alignment: two hands before and after the alignment, respectively at the left and at the right. Each hand is approximated by an ellipsoid. The rotation used apply the axis of one ellipsoid to the axis of the other.}
 		\label{rotation_main}
 		\end{figure}

There are many ways to find the best rotation that aligns two shapes. 
In the case of elongated objects (which was the case in our experiments), one can do the following.
Given two shapes $S_1$ and $S_2$, find the best ellipsoids $\mathbb{E}_1$ and $\mathbb{E}_2$ that approximate the cloud of points defining $S_1$ and $S_2$ respectively, and the unitary matrices $U_1$ and $U_2$ that map the reference axes
to the axes of the ellipsoids (with decreasing lengths). 
The unitary matrices $U_1$ and $U_2$ are uniquely defined if the approximating ellipsoids are triaxial (i.e. the lengths of their principal axes are distinct). Then, we can apply the product matrix $U_2 U_1^{-1}$ on the shape $S_1$ to rotationally align with $S_2$. 
If one encounters a 180 degree flip, apply instead $U_2RU_2^{-1}$, where $R$ is the $180$ degree rotation around the $z$-axis. 
As an example, 
Fig. \ref{rotation_main} shows two hands that have different orientations in space, the corresponding ellipsoids, and the hands after rotation (with a gap to separates them in order to facilitate visualization).
        
\begin{figure}[!ht]
		\centering
 		\includegraphics[width=7cm]{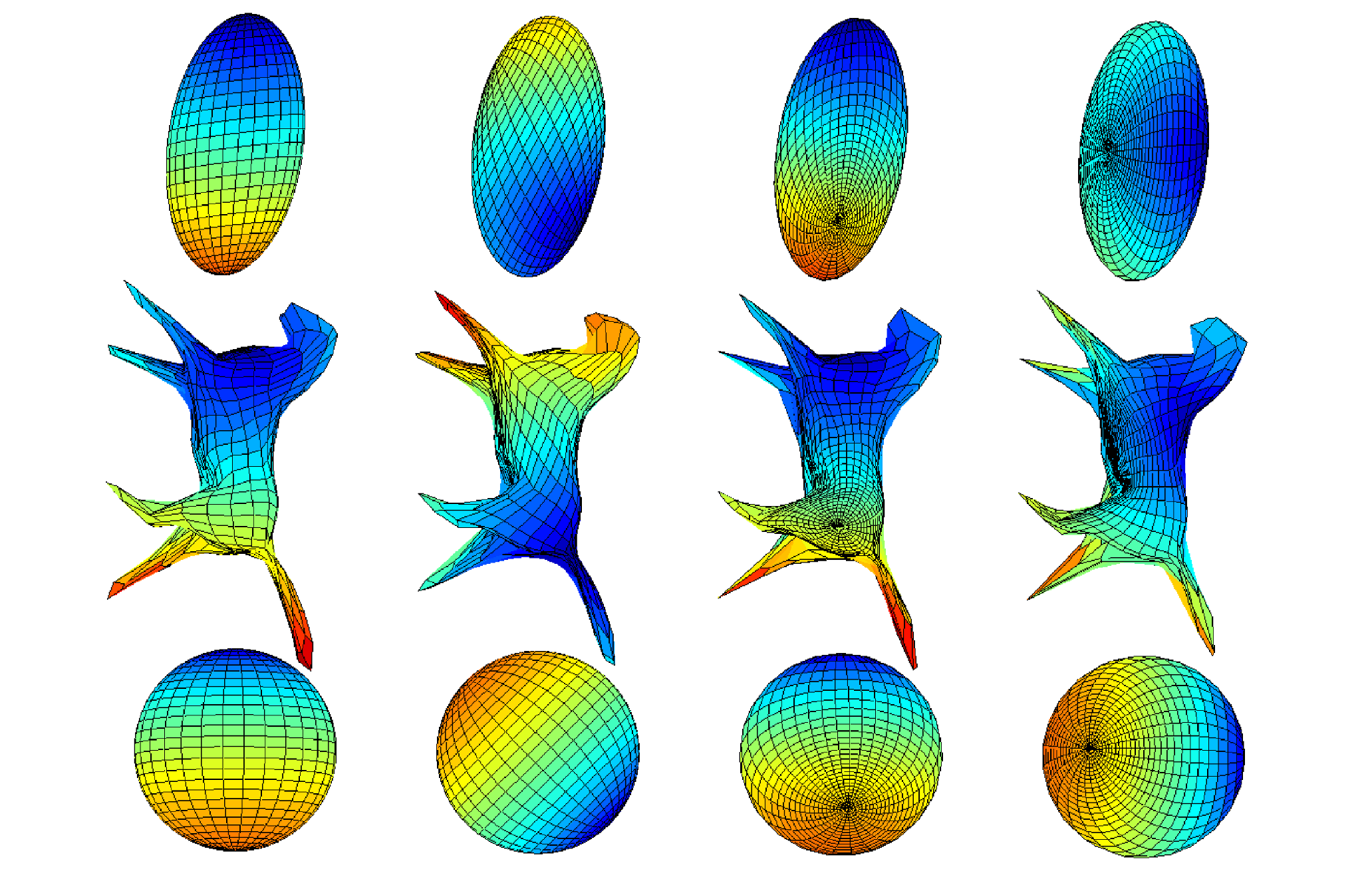}
 		\caption{\scriptsize Robustness of the approximating ellipsoid of a surface with respect to reparameterizations.} 
 		\label{horse_ellipsoid}
 \end{figure}

To find the best ellipsoid that approximates a surface $S$ and the corresponding rotation $U$, one can use a singular value decomposition of $S^{T}S$. However, in the case where the surface is the boundary of a 3D-volume,  
it is more accurate to compute the mean of $S^TS$ over the inscribed volume. 
It also has a more physical meaning since the resulting ellipsoid is equivariant with respect to affine transformations (see 
Section 3 of the Supplementary  Material
where the effect of the rotations of the initial surface on the ellipsoid is illustrated (Fig. 2) and where detailed Algorithms are presented).
 Moreover, the estimation of ellipsoid for an inscribed volume is more stable under reparameterizations. 
 To illustrate this robustness we show in Fig. \ref{horse_ellipsoid} different parameterizations of a horse (middle row) obtained by pre-composing a given parameterization by a diffeomorphism of the sphere (bottom row) and the resulting ellipsoid (top row). 
 The diffeormorphims used in this experiment are (from left to right)
 $\varphi_1 = \textrm{identity}$, 
 $\varphi_2=\textrm{ rotation of} -3\pi/4$ around $x$-axis,
  $\varphi_3 = \textrm{M\"obius transformation that maps } z\in \mathbb{S}^2\simeq\mathbb{C}\cup\{\infty\}$ to $\phi_3(z) =0.4 z + 0.5$,
  $\varphi_4 = $ rotation of $-\pi/2$ around $x$-axis composed with $\varphi_3$.


 In the case where the approximating ellipsoids are not triaxial, one has to use additional information about the surfaces to align them properly (for instance, 
 one can use four points on each surface). This case was not implemented in the present paper.

\subsection{Computations of the energy}
\begin{figure}[!ht]
 		\centering
 	\includegraphics[width=8cm]{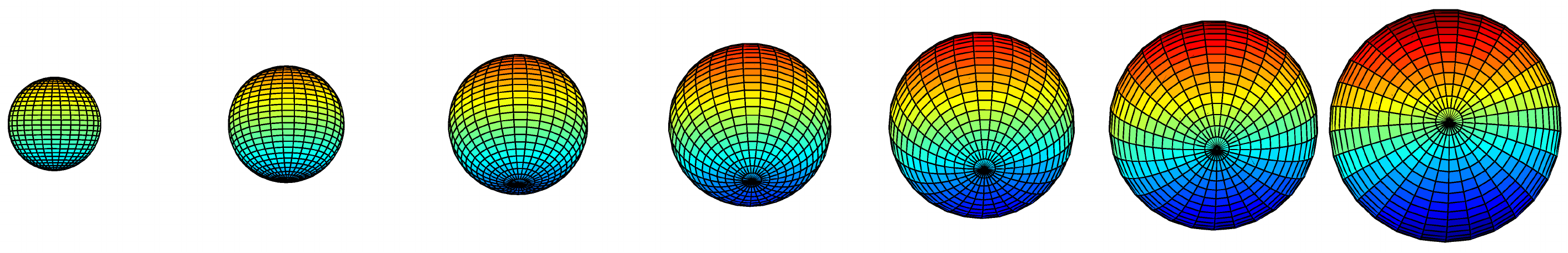}
 		\caption{\scriptsize Path connecting two concentric spheres used for computations in Table \ref{energy1}.
        }
	\label{concentric_spheres}
 		\end{figure}

Let $\Psi:[0,1] \to \mathcal{F}$. The energy of the path $\Psi$ is defined to be:
\begin{align*}
\mathcal{E}(\Psi(t)) & =  \int_{0}^{1} \left\langle\!\!\left\langle \Psi_t^{\perp}, \Psi_t^{\perp} \right\rangle\!\!\right\rangle_{\Psi(t)} dt
= \int_{0}^{1} \left(\! \left( \Psi_t, \Psi_t\right)\!\right)_{\Psi(t)}  dt,
\end{align*}
where $\langle\!\langle \cdot,\cdot \rangle\!\rangle $ is the elastic metric given in Eqn.~\eqref{eqn:elastic-metric2}, $\Psi_t^{\perp} = \left(\Psi_t\cdot n\right) n$ is the normal component of the deformation, and $(\!(\cdot,\cdot)\!)$ is the inner product presented in Eqn.~\eqref{expression_metric}. 
We will present several numerical strategies for approximating this energy and will
compare their computational costs in Table \ref{energy1}. This evaluation uses 
a linear path connecting two concentric spheres of radius $R_1=1$ and $R_2=2.5$,
with constants $a = 1$, $\lambda = 0.125$ and $c = 0$ for defining energy (see Fig. \ref{concentric_spheres}). The theoretical value of the energy
in this case is given by $E_{th} = 32 \pi (a + \lambda)(R_2-R_1)^2$ and measures exclusively the cost of changing the area of the spheres 
(the first and third term of the metric given in Eqn.~\eqref{eqn:elastic-metric2} vanish in this experiment). 
We expect that improvement in accuracy comes at an increased computational cost, 
and this is indeed the case in the results presented in the Table. 
Note that a time-dependent rotation is applied on the path of spheres, but
  the values of the energy is independant of this rotation.

One way to compute the energy of a path $\Psi$ of shapes is to express it using the coefficients of the first fundamental form.
Consider the mapping $\Psi: \mathbb{S}^2 \times \real \to \mathbb{R}^3$ and define
\begin{equation}\label{first_fundamental_form}
\begin{array}{l}
E = \Psi_u\cdot \Psi_u,
F = \Psi_u\cdot \Psi_v,
G= \Psi_v\cdot \Psi_v,
\end{array}
\end{equation}
and their time derivatives 
$$
\begin{array}{l}
\dot{E} = 2 \Psi_{tu}^{\perp} \cdot \Psi_u,
\dot{F} = \Psi_{tu}^{\perp} \cdot \Psi_v + \Psi_u\cdot\Psi_{tv}^{\perp},
\dot{G} = 2\Psi_{tv}^{\perp}\cdot\Psi_v,
\end{array}
$$
as well as the unit normal field $n:= n_f = \frac{f_u\times f_v}{\| f_u\times f_v\|}$ and the vector field $
w=\Psi_{tu}^{\perp}\times \Psi_v+\Psi_u\times \Psi_{tv}^{\perp}.$ Then, the energy of a path $\Psi$ decomposes into the sum of four terms:
$\mathcal{E}(\Psi(t)) = \mathcal{E}_1+\mathcal{E}_2+\mathcal{E}_3+\mathcal{E}_4$, where 
$$
\begin{array}{ll}
\!\!\!\mathcal{E}_1 \!=\!&\!\!\! a \int_{0}^{1}\!\!\! \int_{\mathbb{S}^2} (EG-F^2)^{-3/2} B\,d\!u\,d\!v\,dt
\\
&
\begin{array}{rl}
 \textrm{ with } B = &\!\!\! G^2 \dot{E}^2 \!+\! 2(EG\!+\!F^2)\dot{F}^2 \!+\! E^2\dot{G}^2 \\
 &\!\!\!-\! 4 FG\dot{E}\dot{F} \!+\! 2F^2\dot{E}\dot{G} \!-\! 4 E F\dot{F}\dot{G}\ ,\end{array}
\\
\!\!\!\mathcal{E}_2 \!=\!&\!\!\!\left(\frac{\lambda}{2}\!+\!\frac{c}{4}\right)\!\int_{0}^{1}\!\!\!\int_{\mathbb{S}^2} 
\!(EG\!-\!F^2)^{{-\frac{3}{2}}}\!(G\dot{E}\!-\!2F\dot{F}\!+\!E\dot{G})^{{2}}\!d\!u\,d\!v\,dt\ ,
\\
\!\!\!\mathcal{E}_3 \!=\! &\!\!\!-c \int_{0}^{1} \int_{\mathbb{S}^2}  (G\dot{E}-2F\dot{F}+E\dot{G})  (n\cdot w)\,d\!u\,d\!v\,dt\ ,
\\
\!\!\!\mathcal{E}_4 \!=\!&\!\!\!  c \int_{0}^{1} \int_{\mathbb{S}^2}(EG-F^2)^{-\frac{1}{2}}  (w\cdot w)\,d\!u\,d\!v\,dt \ .

\end{array}
$$
In the implementation of these formulas, we can reach singularities on the boundary of the integration domain, which we can ignore. 
In the example involving two concentric spheres, 
the total energy computed by this method is labelled $E_{I\&II}$ in Table \ref{energy1}.

Another way to compute the energy is based on Eqn.~\eqref{expression_metric} that expresses the 
elastic metric in terms of principal curvatures. 
In terms of the coefficients of the first fundamental form given in Eqn.~\eqref{first_fundamental_form} and of the second fundamental given by
$$
\begin{array}{l}
e = \Psi_{uu}\cdot n = - \Psi_u\cdot n_u,\\ f = \Psi_{uv}\cdot n = -\Psi_u\cdot n_v = -\Psi_v\cdot n_u,\\ g= \Psi_{vv}\cdot n = - \Psi_v\cdot n_v,
\end{array}
$$
the Gauss curvature $K$ and the mean curvature $H$ have the following expressions
$$
K = \frac{eg-f^2}{EG-F^2}, \quad H = \frac{1}{2}\frac{eG + gE - 2 fF}{EG-F^2},
$$
and the principal curvatures are given by
$$\kappa_1 = H + \sqrt{H^2 - K}, \quad \kappa_2 = H-  \sqrt{H^2 - K}.$$ 
Again, in the implementation of these formulas, we can get singularities for curvatures on the boundary, but we can ignore them in computing 
the integral given in
Eqn.~\eqref{expression_metric}.
This corresponds to removing a small disc on the parameterized surface around the images of the north and south poles. 

In the example of the two concentric spheres, the theoretical values of $\kappa_1$ and $\kappa_2$ is the constant function equal to $1/R$ where $R = R_1 + t(R_2 - R_1)$ is the radius of the sphere along the path interpolating linearly the sphere of radius $R_1 = 1$ to the sphere of radius $R_2 = 2.5$.
The total energy computed by this method is labelled $E_{k_1k_2}$ in Table \ref{energy1}. 

\begin{figure}[!ht]
		\includegraphics[width=7cm]{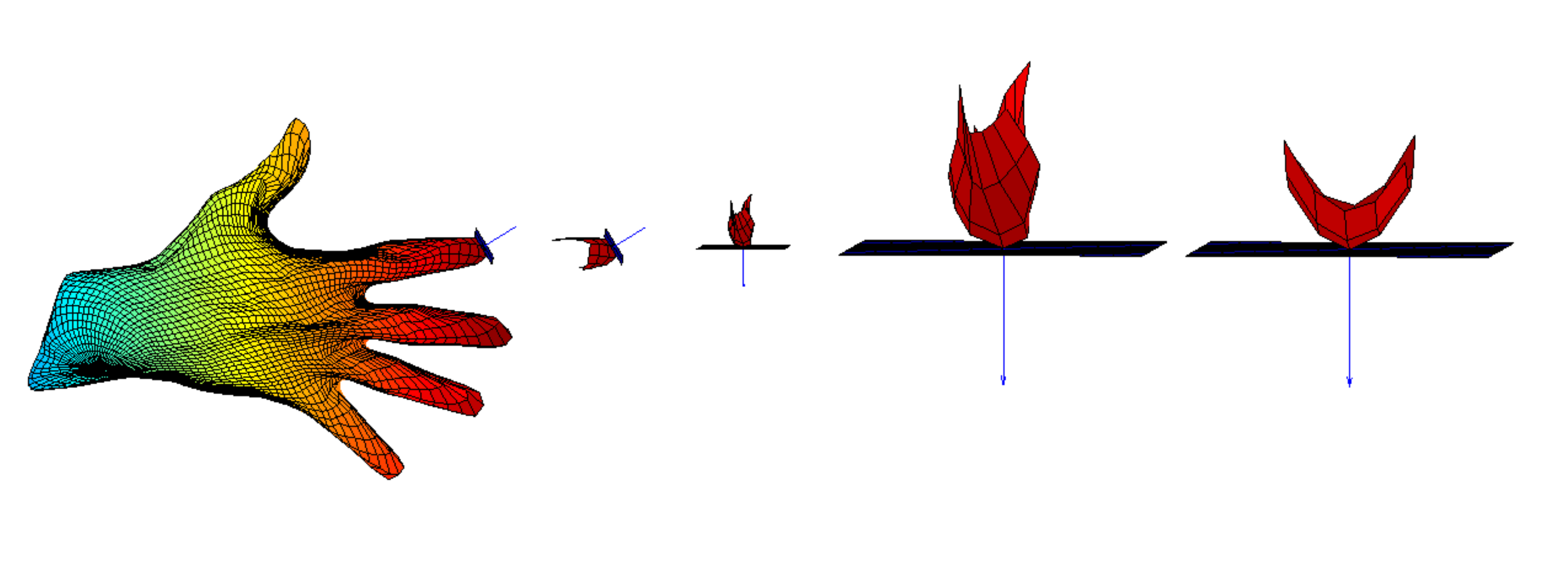}
 		\caption{\scriptsize From left to right: A hand with the tangent plane and normal at the tip of the index finger; 3-neighborhood of the tip of the index finger; tip of the index finger after rotation; a closeup; approximating second order polynomial.}
 		\label{courbure_main}
 		\end{figure}

To improve the computation of the curvatures and therefore also of the energy, we can use 
polynomial approximations of the surfaces. This procedure, leading to the computation of the principal curvatures, is illustrated in Fig.~\ref{courbure_main}. To compute the principal curvatures at a given point of a surface, e.g. at the tip of the index finger of the hand depicted in Fig.~\ref{courbure_main}, we first compute the normal at this point by averaging the normals of the facets having this point as vertex. A tangent plane is then defined as the plane orthogonal to the normal passing through the point under consideration. A neighborhood of the point is isolated from the surface (we use a 3-neighborhood, see second drawing in Fig.~\ref{courbure_main}). We then apply a rigid transformation to center the point at the origin and to align the tangent plan with the xy-plane (see third drawing, and a closeup in 
the fourth drawing). After that, we use Algorithm~5
given in Section~5 of the Supplementary Material
to compute the second order polynomial $P(x,y) = a_1x^2+a_2y^2+a_3xy+a_4x+a_5y+a_6$, which minimizes the sum $\sum_i(z_i - P(x_i, y_i))^2$ over the points of the centered and rotated neighborhood. Then, the Gauss curvature at that point is given by $K = 4a_1a_2-a_3^2$, the mean curvature by $H = a_1+a_2$, and the principal curvatures by $\kappa_1 = a_1+a_2+\sqrt((a_1-a_2)^2+a_3^2)$ and $\kappa_2 = a_1+a_2-\sqrt((a_1-a_2)^2+a_3^2) $. In the example of the two concentric spheres, 
 the total energy computed using the principal curvatures obtained by this method is labelled $E_{P}$ in Table \ref{energy1}.



\begin{figure}[!ht]
\includegraphics[width=7cm]{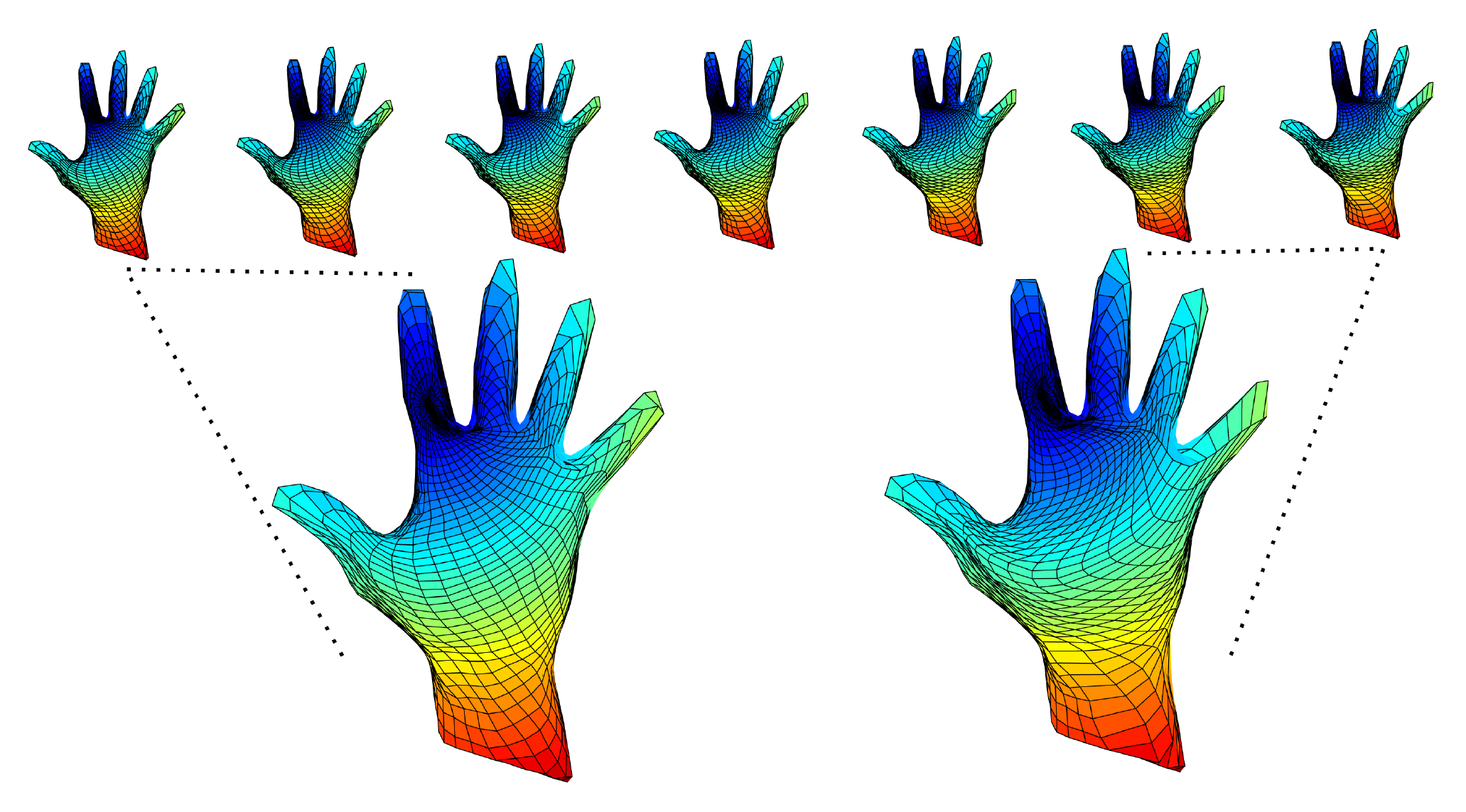}
 		\caption{\scriptsize 
		A path of zero energy connecting a hand and the same hand with another parameterization.  
		}
 		\label{Picture7_hand_2_hand}
 		\end{figure}

In order to show that the energy function of a path of shapes is independent of the way the objects 
are parameterized, we replace the integration over the domain of parameterization by the integration over the triangulated surfaces. This means that we approximate
the area elements of the surfaces by the area of triangles whose vertices are given by the parameterization. In this way, the parameterization of surfaces is only used to define the surfaces, but plays no role at all in the computation of the energy function. In the example of the two concentric spheres, 
the total energy computed by this method is labelled $E_{\Delta}$ and given in Table \ref{energy1}. 
In Fig.~\ref{Picture7_hand_2_hand}, a path connecting a hand to the same hand, but with a different parameterization, is shown. The energy  of this path, computed with the constants $a=1$, $\lambda = c = 0.125$, reads $E_\Delta = 0.4824$, hence is close to $0$.  
Now returning to Fig. \ref{fig_illustration}, the energy of the lower path from a horse to a cat computed 
with the same constants is $E_\Delta = 227.4049$, its length is $L[\Psi] = 14.9099$, whereas the upper path 
(obtained from the lower path by applying a different reparameterization at each time step) has an energy equal to $E_\Delta = 225.5249$ 
and a length of $L[\Psi] = 14.8802$.
Note that in this example, 
the colors refer to the Euclidean distance to the point on the surface corresponding to the image of the north pole of the sphere 
(cold colors for small distances versus hot colors for large distances). In particular, the north and south poles do not correspond in these two paths.
\begin{table}
\begin{center}
\begin{tabular}{|c|c|}
\hline 
Energy, $10^4$ points per object & Elapsed time for $10^4$ points \tabularnewline
\hline
\hline 
$E_{I\&II}  =   246.2854$ &  0.221726 seconds  
\tabularnewline
\hline 
$E_{k_1k_2}= 249.1969$ & 0.862376 seconds
\tabularnewline
\hline 
$E_P = 255.8288
$ & 1.238354
seconds
\tabularnewline
\hline 
$E_\Delta = 255.9043$ & 9.738431 seconds
\tabularnewline
\hline
\hline
  Energy  for $4\times 10^4$ points & Elapsed time for $4 \times10^4$ points\tabularnewline
\hline
\hline 
 $E_{I\&II}  =   249.1503$ &  0.978828 seconds
 \tabularnewline
\hline
$E_{k_1k_2}= 251.8494$ & 3.45599 seconds
 \tabularnewline
\hline
$E_P = 254.7646
$ & 4.906798
seconds 
\tabularnewline
\hline
$E_\Delta = 254.7832$ & 39.011899  seconds
 \tabularnewline
\hline
\end{tabular}
\end{center}
\caption{\scriptsize Computation of the energy of a path connecting two concentric spheres (Fig. \ref{concentric_spheres}) using different methods, and time needed for the computations. The theoretical value of the energy is $E_{th} = 254.4690$. Here $R_1 = 1$, $R_2 = 2.5$, $\lambda = 0.125$ and $c = 0$.}
\label{energy1}
\end{table}
\subsection{Orthonormal Basis of Deformations}
In this section, we define bases for representing perturbations of a path of surfaces. These basis elements form possible directions for use in path-straightening in
Section \ref{Path-straightening}.
The first basis 
we used is a variation of the one given in \cite{KurtekPAMI2012}. We start with a basis $\mathcal{B}_1 = \{Y_{l}^m, 1\leq l \leq N, -l\leq m\leq l\}$ of spherical harmonics of degree less than $N$, available in Matlab as function SPHARM (see \cite{courant_hilbert} for more information on spherical harmonics). We make three copies of this basis of $\mathbb{R}$-valued functions in order to obtain a basis $\mathcal{B}_2$ of the space $L^2(\mathbb{S}^2, \mathbb{R}^3)$ of $\mathbb{R}^3$-valued functions. Similar to Xie et al. \cite{XieICCV2013}, we demonstrate 
reconstruction of some surfaces using the resulting basis, as the degree of the spherical harmonics grows, in the Supplementary 
Material (Fig. 5).


Next, we want to construct a basis of perturbations of a path connecting two parameterized surfaces $f_1$ and $f_2$. In order to apply the path-straightening method as described in Section \ref{Path-straightening}, we want the perturbations to vanish at $t = 0$ and $t=1$ so that $f_1$ and $f_2$ remain fixed. Therefore, 
we want a basis of $L^2(\mathbb{S}^2\times [0,1], \mathbb{R}^3)$ with elements that have this property. To ensure this, each element of $\mathcal{B}_2$ is multiplied by a basis element of $L^2([0, 1],\mathbb{R})$ of the form $P_j(t) = \frac{1}{4}\sin(\pi jt)$, $1\leq j\leq J$. 
Unfortunately a major limitation of the resulting $L^2$ basis is that slowly- and  rapidly-oscillating harmonics have comparable amplitudes. In the implementation of the path-straightening method, this implies that the updated path can go out of the open set of  immersions. 

One possible way to counter this effect is to orthonormalize the $L^2$-basis with respect to an $H^1$-type scalar product 
(i.e. that measures also the variation of the derivatives). For this kind of scalar product, an orthonormal basis consists of 
functions which have controlled derivatives (hence can not oscillate to much). This approach was also used in \cite{KurtekPAMI2012} 
where the $L^2$-basis is orthonormalized   using the Gram-Schmidt procedure with respect to the following scalar product
$$
\begin{array}{ll}
(B^1, B^2) =&\!\!\!\! \int_{0}^1\int_{\mathbb{S}^2}\left(B^1\cdot B^2 + B^1_t\cdot B^2_t+ B^1_u\cdot B^2_u\right. \\&\!\!\!\! \left.+ B^1_v\cdot B^2_v + B^1_{t,u} \cdot B^2_{t,u}+ B^1_{t,v} \cdot B^2_{t,v}\right)ds\,dt.
\end{array}
$$
However, when increasing the degree of spherical harmonics, the computational cost of generation of an orthonormal basis using this scalar product is very high. 
Therefore, we first orthonormalize the basis $\mathcal{B}_2$ with respect to the following inner product 
\begin{equation}\label{H1}
\begin{array}{ll}
(B^1, B^2) =& \!\!\!\!\int_{\mathbb{S}^2}\left(B^1\!\cdot\!B^2 + B^1_u\!\cdot\!B^2_u+ B^1_v\!\cdot\!B^2_v \right)ds,
\end{array}
\end{equation}
and then we multiply the resulting basis by the time-dependant components $P_j(t) = \frac{1}{4}\sin(\pi jt)$,  $1\leq j\leq J$.
The advantage of this method is that the Gram-Schmidt procedure is applied to matrices of lower dimensions (without the time dimension) 
and on a smaller number of elements (by a factor $J$). The spatial oscillations of the resulting basis elements are well controlled by 
the presence of the spatial derivatives $B_u$ and $B_v$ in the inner product given in Eqn.~\eqref{H1}.

\subsection{Path-straightening method}\label{Path-straightening}
The path-straightening method is used to find critical points of the energy functional.
Starting with an arbitrary path, the method consists of iteratively deforming (or ``straightening'') the path in the opposite direction of the gradient, until the path converges
to a geodesic. 
The gradient of the path energy is approximated using 
a basis $\mathcal{B}$ of possible perturbations of a path of surfaces $\Psi$, as constructed in the previous section. 
We first compute the directional derivatives  
$\nabla \mathcal{E}_{\Psi}(b)=\frac{d}{d\mathcal{\epsilon}}(\mathcal{E}(\Psi + \epsilon b))|_{\epsilon=0}
$
where $b$ ranges over $\mathcal{B}$. 
This is done by fixing a small $\epsilon_1$ and approximating the directional derivative by
$
\nabla \mathcal{E}_{\Psi}(b)\simeq (\mathcal{E}(\Psi + \epsilon_1 b)- \mathcal{E}(\Psi)){\epsilon_1}^{-1}.
$
Using the finite orthonormal basis $\mathcal{B}$, we obtain a numerical approximation of the gradient:
$
\nabla \mathcal{E}_{\Psi} = \sum_{b\in\mathcal{B}}\nabla \mathcal{E}_{\Psi}(b)~b.
$
In particular, the norm of the gradient is approximately given 
by $\|\nabla \mathcal{E}_{\Psi}\|^2 = \sum_{b\in\mathcal{B}}\nabla \mathcal{E}_{\Psi}(b)^2$.
The update of the path is done by replacing $\Psi$ by $\Psi - \epsilon_2\nabla \mathcal{E}_{\Psi}$, 
where $\epsilon_2$ is a small parameter that has to be ajusted empirically. 
The method is detailed in Algorithm \ref{alg:psi_point} below.


\begin{algorithm}[ht]\label{algo}
\begin{scriptsize}
\KwIn{\begin{enumerate}\item A path  $\Psi$ between two parameterized surfaces ${f}_1$ and ${f}_2$,\item a basis of perturbation $\mathcal{B}$.\end{enumerate}}
\KwOut{\begin{enumerate}\item The minimal energy needed to deform $f_1$ into $f_2$ given by the value of the cost function $E$, \item the geodesic path between $f_1$ and $f_2$. \end{enumerate}}

Set $\|\nabla{E}\|^2 = 1$.

\While{$\|\nabla{E}\|^2 >  10^{-3}$}
{
\textbf{2-} Compute  the energy $E$ of the path $\Psi$ according to Eqn.\eqref{eqn:elastic-metric2} or Eqn.~\eqref{expression_metric}.

\textbf{3-} Set $\Psi_{\textrm{upd}}=0$ and $\| \nabla E\|^2 = 0$.

\For{$i\leftarrow 1$ \KwTo $size(\mathcal{B})$}
{
\textbf{4-} Add a perturbation to the current path $\Psi$: define $\Psi(i) = \Psi + \epsilon_1\, \mathcal{B}(i)$, 
where $\mathcal{B}(i)$ is the element of the perturbation basis $\mathcal{B}$ of index $i$ and $\epsilon_1>0$ is small.\\
\textbf{5-} Compute  the energy $E(i)$ of the perturbed path~$\Psi(i)$.\\
\textbf{6-} Compute the gradient of energy $\nabla{E}(i)$   in the direction $\mathcal{B}(i)$ using 
the approximation $\nabla{E}(i)\sim\frac{E(i)-E}{\epsilon_1}$.

\textbf{7-} Compute the updating path: $\Psi_{\textrm{upd}}\leftarrow \Psi_{\textrm{upd}}+\nabla{E}(i)\cdot \mathcal{B}(i)$.

\textbf{8-} Compute the squarred norm of the gradient of energy at path $\Psi$:
$
\|\nabla{E}\|^2 \leftarrow \|\nabla{E}\|^2 + (\nabla{E}(i))^2.
$
}
\textbf{10-} Update the path: $\Psi = \Psi - \epsilon_2  \Psi_{\textrm{upd}}$

}


\end{scriptsize}
\caption{Path-straightening method.}
\label{alg:psi_point}
\end{algorithm}
\vspace{-0.5cm}
\section{Experimental results}

The 3D realistic models used in our experiments are part of the TOSCA \cite{B2008} 
dataset. 
Their spherical parameterizations were initially implemented
in \cite{KSKH13}.


\subsection{Examples of geodesics obtained by path-straightening}

\begin{figure}[!ht]
 		\centering
 		\includegraphics[width=8cm]{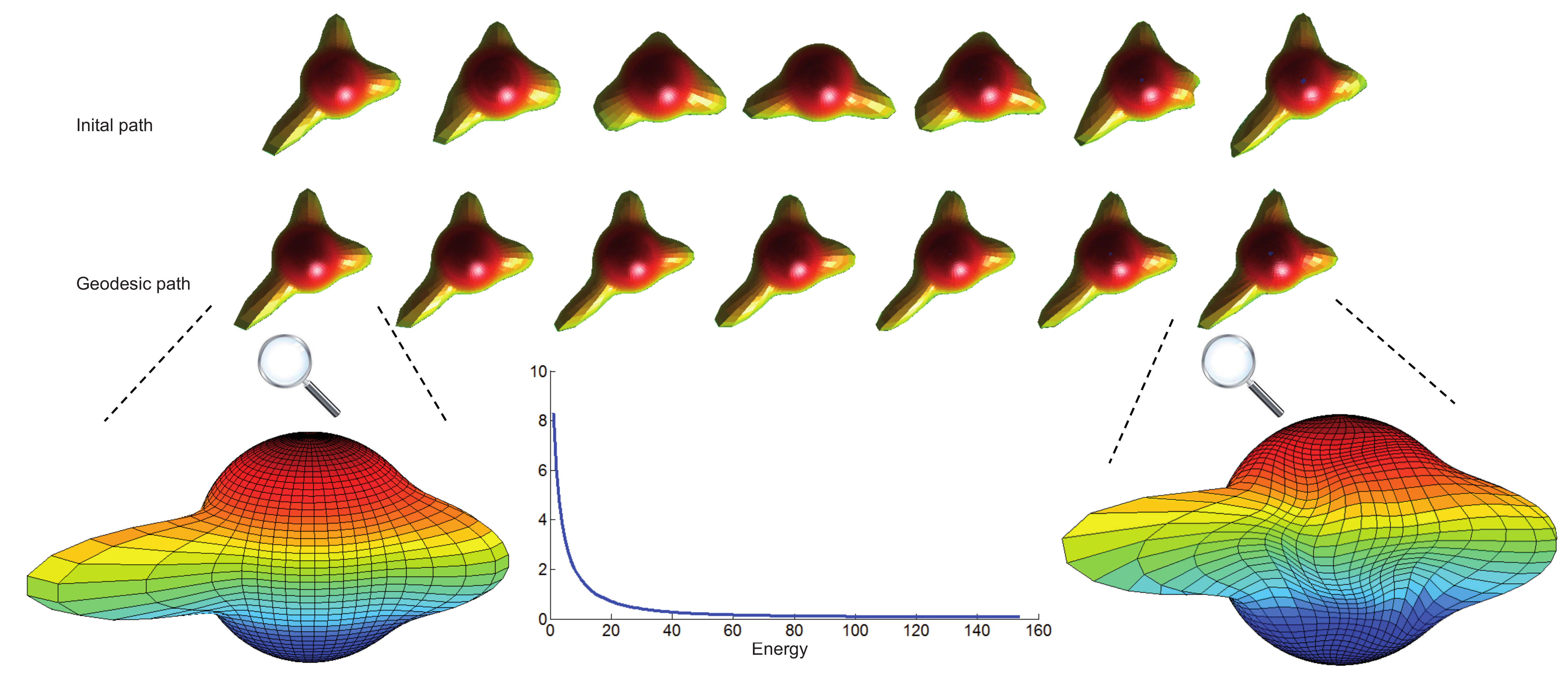}
 		\caption{\scriptsize Illustration of initial path (upper row) and geodesic path in shape space (middle row). The energy is reported in the buttom row. The surfaces at the end points of the path have different parameterizations. 
        }
 		\label{illustrer_gauge}
\end{figure}

First we apply the path-straightening method to the case where 
the surfaces at the extremes of the initial path have the same shape, but different parameterizations. More precisely, 
we consider the special case where $\Psi_0(0)=f_1$,  $\Psi_0(1)=f_1\circ \gamma$ for some diffeomorphism $\gamma$ 
and where we initialize the path with 
piecewise linear interpolation to a different surface $f_3$ 
in the middle of the path, i.e.  $\Psi_0(\frac{t}{2})=f_3$.	
This situation is illustrated in Fig. \ref{illustrer_gauge}. The proposed gauge-invariant approach is expected to reach a path with
constant shape as a geodesic, despite the different shapes appearing in the initial path and  the different parameterization of shapes at the 
end points of the path
(to emphasize the differences in parameterization, zoom-ins of these surfaces are also shown).
Once we have the geodesic path $\Psi$ between the given surfaces, the distance in the shape space between $f_1$ and $f_1 \circ \gamma$, 
$d_{\Psi}(f_1, f_2)$, is the length of $\Psi$ as specified in Eqn.~\eqref{length}.
As expected,  the resulting geodesic path, shown in Fig.~\ref{illustrer_gauge},  is constant with the same shape as the either end, 
and with $d_{\Psi}(f_1, f_2)=0$.  Using 
path-straightening, we obtain a $99.28\%$ decrease in the energy function from the initial path to the final path. 

\begin{figure*}[!ht]
		\centering
            \includegraphics[width=15cm]{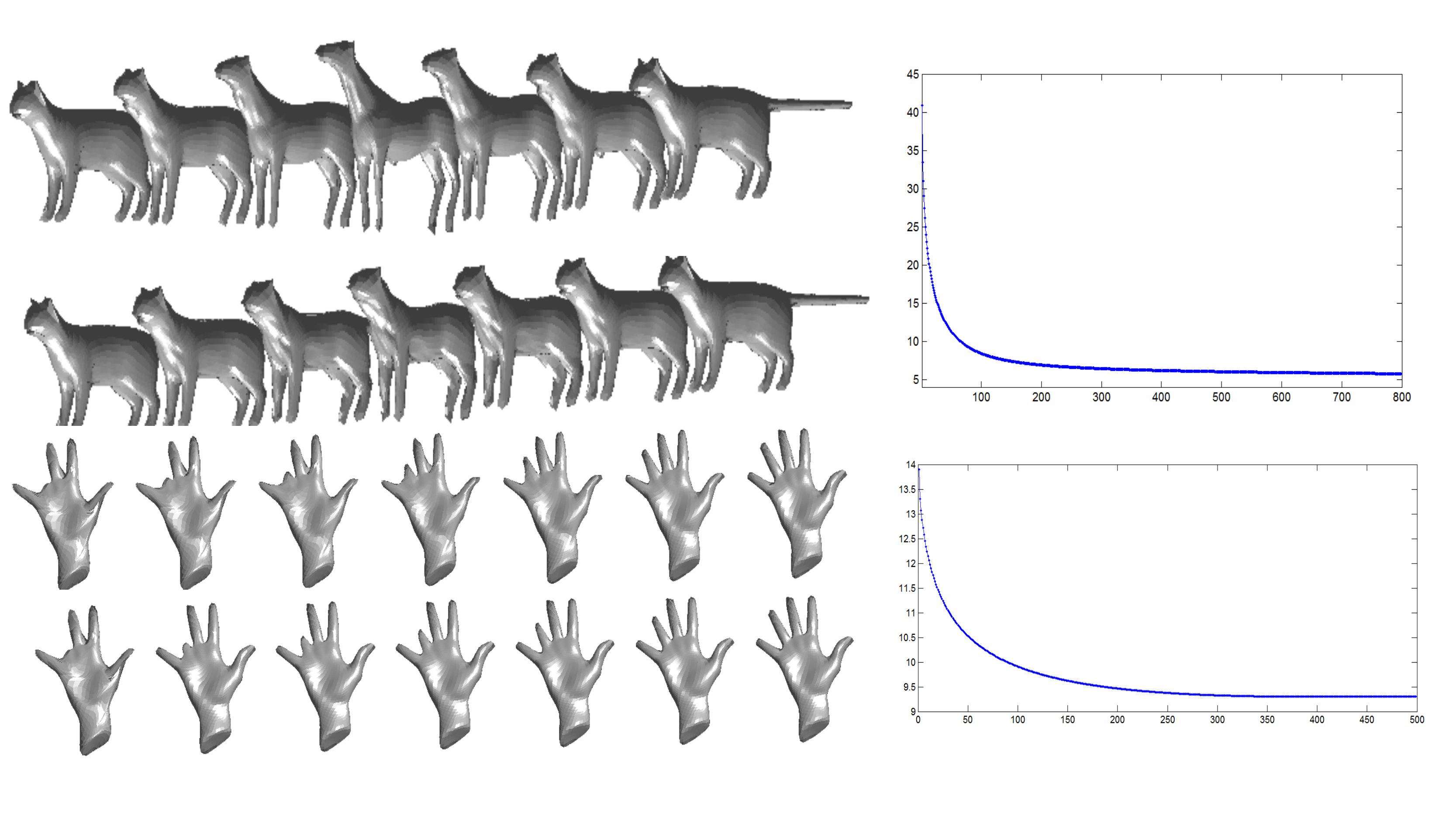}
 		\caption{\scriptsize The top row shows an initial path formed by linear interpolation between a cat to a horse and back to the cat. The second row illustrates the geodesic obtained after 800 iterations of path-straightening. The corresponding evolution of the energy is shown on the right. 
		Similarly, the third row shows a linear path between two hands with bad correspondence and the last row shows the final geodesic, with the
		corresponding energy is shown on the right.}
 		\label{geodesic_nv}
 \end{figure*}

In Fig. \ref{geodesic_nv} we consider more challenging shapes. The top-two rows display the case where we have $\Psi_0(0)=f_1$,  $\Psi_0(1)=f_1$ (a cat) 
and where we initialize the path with piecewise linear interpolation to a horse in the middle of the path. 
The upper row shows the initial path and the second row the geodesic path. We can see that the geodesic path has a constant
shape throughout, as expected. We also plot the evolution of the 
path energy on the right during path-straightening. We can see that
the energy decreases until it reaches a relatively small value; the theoretical minimum is, of course, zero for a contant path. In the last two rows of Fig. \ref{geodesic_nv}, we consider the case of two hands.
We initialize the path with linear interpolation (third row in Fig. \ref{geodesic_nv}), and the resulting path is shown in 
the last rows of Fig. \ref{geodesic_nv}. 
The energy evolution is shown on the right and we can see the energy decreasing until it reaches a constant value; 
thus, the final path is a geodesic. It can be seen that the deformation along the geodesic path is more natural than the original path.

\subsection{Classification of 3D shapes}
 As mentioned earlier, the geodesic paths provide us with tools for comparing, and deforming parameterized surfaces. We suggest a comparison of shapes of 3D objects using geodesic distances between their boundary surfaces in the shape space. This section presents a specific application to illustrate that idea.
 In this section, we study several shapes belonging to four classes: horses, hands, cats and centaurs. 

\begin{figure}[!ht]
  		\centering
 \includegraphics[width=9cm]{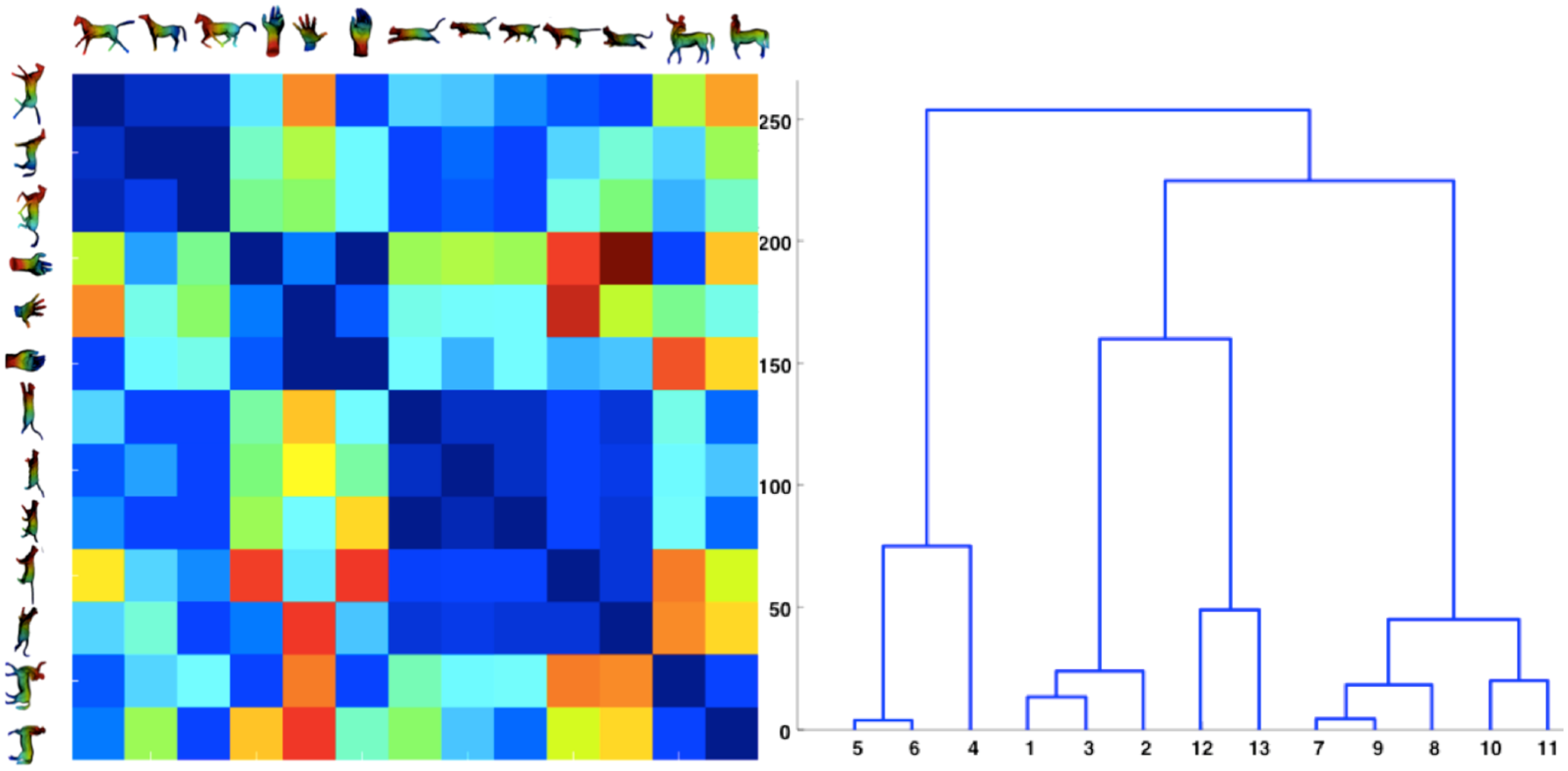}
  		\caption{\scriptsize Classification performance; left: the distance matrix. right: the dendrogram.}
  		\label{fig_classif}
  		\end{figure}

We begin by computing the pairwise geodesic distances between corresponding 3D surfaces. The distance matrix and the classification dendrogram are shown in Fig. \ref{fig_classif}. In the distance matrix, we can easily distinguish four classes corresponding to four blue boxes. Actually the cold colors in the illustrated matrix correspond to small values of distances versus hot colors that correspond to greater distances. 
The clustering obtained using the {\it dendrogram} (command in matlab) can be interpreted by slicing the top of the dendrogram by a horizontal line to split the shapes into the desired number of classes, and then sliding the horizontal line to the bottom in order to refine the classification. The coarsest classification results by slicing the dendrogram into two classes (by a horizontal line close to the top), the shapes 4, 5 and 6 (the hands) forms a first class and the remaining (horses, cats and centaurs) are grouped together as a second class. The next level in classification distinguishes the shapes 1, 2, and 3 (the horses) and 12, 13 (the centaurs) from the shapes 7, 8, 9, 10, 11 (the cats). The finest level separates the horses and the centaurs in different classes and results in four classes.
Thus, we argue that the proposed framework provides a powerful tool for shape classification.


\subsection{The effect of number of basis elements}

\begin{figure}[!ht]
 		\centering
\includegraphics[width=8cm]{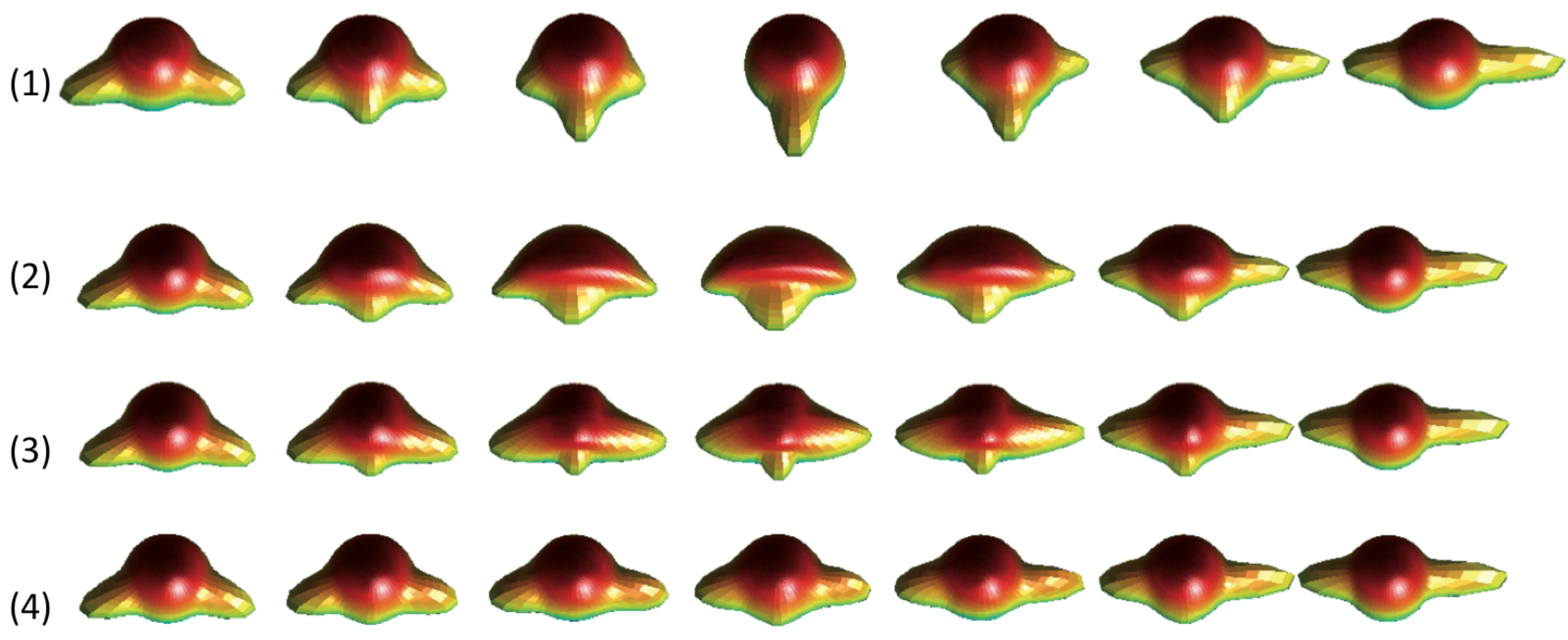}
\includegraphics[width=8cm]{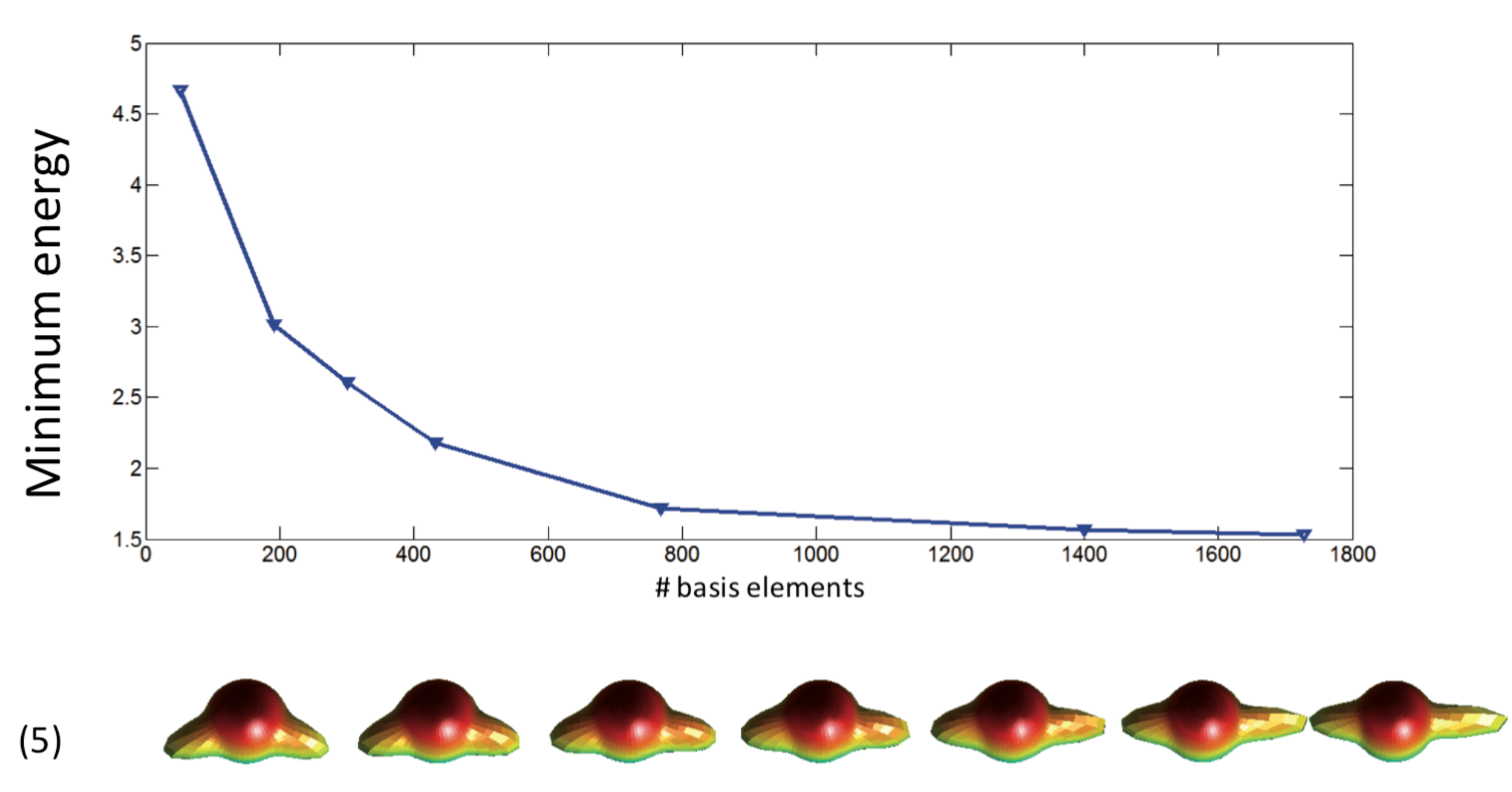}
 		\caption{\scriptsize The effect of the number of basis elements, (1) initial path, (2) geodesic path using 52 basis elements, (3) geodesic path using 432 basis elements, (4) geodesic path using 1728 basis elements, (5) geodesic path using 1728 basis elements after linear interpolation initialization.}
 		\label{fig_basis}
\end{figure}

In this section, we study the effect of the number of basis elements, used
in path-straightening, on the resulting geodesic path. Given two parameterized surfaces $f_1$ and $f_2$, 
we again initialize the path with the linear interpolation to a different surface $f_3$ in the middle of the path. 
This initial path is shown in the upper row of Fig \ref{fig_basis}. Then, we compute the geodesic path using different 
number of basis elements. We show the geodesic paths that use $52$, $432$ and $1728$ basis elements, respectively. We can see that the 
larger the number of basis elements, the better the final result is. We also provide the trade-off between the number of 
basis elements and the minimum energy value obtained. The trade-off confirms our assertion. At the bottom of the figure, 
we show the geodesic path obtained when the path-straightening Algorithm is initialized with the linear interpolation 
between $f_1$ and $f_2$. This path is also calculated using the number of basis elements corresponding to the lowest energy. 
This path can be seen as ground truth to visually interpret the previous geodesics (with more complicated initial conditions and fewer basis elements). 

%
%
%

\section{Conclusion}

In this paper we have proposed a novel Riemannian framework for computing geodesic paths between shapes of parameterized surfaces.
These geodesics are invariant to rigid motion, scaling and most importantly reparameterization of individual surfaces. The novelty lies in defining a Riemannian metric directly on the quotient (shape) space, rather
than inheriting it from pre-shape space, and in using it to formulate a path energy that measures only the normal components of velocities along the path. The geodesic
computation is based on a path-straightening technique that iteratively corrects paths between surfaces until geodesics are achieved. 
We have presented some examples of geodesics between surfaces in shape spaces and utilized the distances between surfaces for classification of some 3D shapes. 
However, the computational costs of our programs are deemed high and convergence should be accelerated in order to be able to apply this framework in realistic 
practical scenarios such as, for instance, human body action recognition.

%
\ifCLASSOPTIONcompsoc
 \section*{Acknowledgments}
\else
 \section*{Acknowledgments}
\fi

The authors would like to thank the anonymous reviewers for valuable comments that helped improve the paper significantly. 
The spherical parameterizations of the surfaces used in this paper were implemented by H. Laga. 
The authors would like to thank S.~Kurtek for his help, and J.-C. Alvarez-Paiva and M. Bauer for  fruitful discussions. 
This work was supported in part by the Labex CEMPI  (ANR-11-LABX-0007-01), the Equipex IrDIVE (ANR-11-EQPX-23 “IrDIVE”) and was made possible 
by a visit of ABT
to Telecom-Lille CRIStAL financed by the CNRS. During the reviewing process, ABT enjoyed excellent working conditions 
at the Pauli Institute, Vienna, Austria. Both ABT and AS benefitted from the program on Shape Analysis held 
in 2015 at the ESI, Vienna, Austria.

%

\ifCLASSOPTIONcaptionsoff
  \newpage
\fi

\bibliographystyle{IEEEtran}
\bibliography{ShapeBiblio}
%

\newpage

\parpic{\includegraphics[width=0.9in,clip,keepaspectratio]{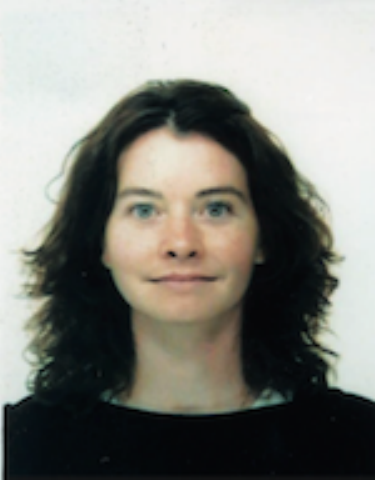}}
\noindent \begin{scriptsize}{{\bf Alice Barbara Tumpach}  is an Associate Professor in Mathematics (University Lille 1, France) and member of the Laboratoire Painlev\'e (Lille 1/CNRS UMR 8524), since 2007. She received a Ph.D degree in Mathematics in 2005 at the Ecole Polytechnique, Palaiseau, France. She spent two years at the Ecole Polytechnique F\'ed\'erale de Lausanne as a Post-Doc. Her research interests lie in the area of infinite-dimensional Geometry, Lie Groups and Functional Analysis.}\end{scriptsize} 

\parpic{\includegraphics[width=0.9in,clip,keepaspectratio]{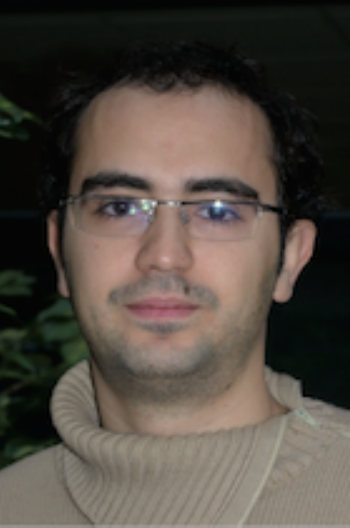}}
\noindent \begin{scriptsize}{{\bf Hassen Drira} is an assistant Professor of Computer Science at T\'el\'ecom Lille and member of the Laboratoire CRIStAL (UMR CNRS 9189) since September 2012. He obtained his Ph.D degree in Computer Science in 2011 from University of Lille 1 (France). His research interests are mainly focused on pattern recognition, statistical shape analysis. He has published several refereed journals and conference articles in these areas.}\end{scriptsize}

\parpic{\includegraphics[width=0.9in,clip,keepaspectratio]{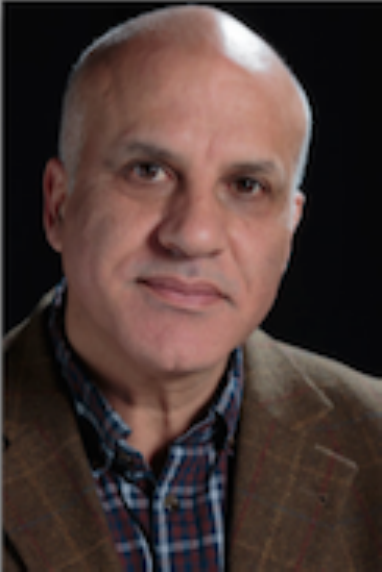}}
\noindent \begin{scriptsize}{{\bf Mohamed Daoudi} is a Professor of Computer Science at Telecom Lille  and member of the Laboratoire CRIStAL (UMR CNRS 9189). He received his Ph.D. degree in Computer Engineering from the University of Lille 1 (France) in 1993 and Habilitation à Diriger des Recherches from the University of Littoral (France) in 2000. His research interests include pattern recognition, shape analysis and computer vision.  He is  a Fellow member of the IAPR.}\end{scriptsize}

\parpic{\includegraphics[width=0.9in,clip,keepaspectratio]{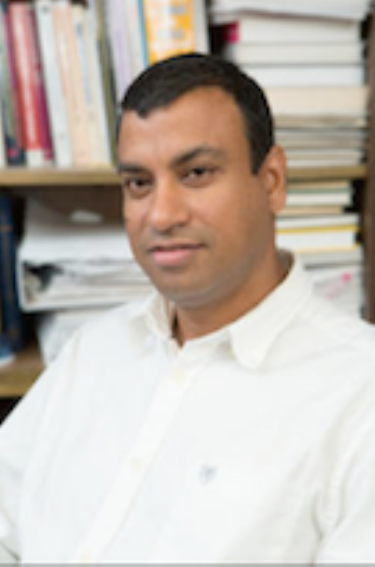}}
\noindent \begin{scriptsize}{{\bf Anuj Srivastava} is a Professor of Statistics at the Florida State
University in Tallahassee, FL. He obtained his MS and PhD degrees in
Electrical Engineering from the Washington University in St. Louis in 1993
and 1996, respectively. After spending the year 1996-97 at the Brown
University as a visiting researcher, he joined FSU as an Assistant
Professor in 1997. His research is focused on pattern theoretic approaches to problems
in image analysis, computer vision, and signal processing. }\end{scriptsize} \\

\newpage

\appendices
\begin{center}
\Large{SUPPLEMENTARY MATERIAL}
\end{center}
\section*{The Subbundle of normal vector fields is $\Gamma$-invariant}
In this section we provide the proof of the fact that the subbundle of normal vector fields is a $\Gamma$-invariant complement to the subbundle of tangent vector fields.

\begin{proposition}\label{Nor_Gamma_invariant}
Denote by $Nor$ the subbundle of the tangent bundle $T\mathcal{F}$ consisting of \textit{normal vector fields} which is the space of vector fields such that the cross product with the normal $n_f~:\mathbb{S}^2\rightarrow\mathbb{R}^3$ to the shape $f(\mathbb{S}^2)$ vanishes~:
$$
Nor(f) = \{\delta f~:\mathbb{S}^2\rightarrow \mathbb{R}^3, \textrm{~such~that~} \delta f\times n_f = 0\}.
$$ 
Any tangent vector $\delta f\in T_f\mathcal{F}$ admits a unique decomposition 
$$
\delta f = \delta f^{T} + \delta f^{\perp}
$$
into its tangential part $\delta f^{T}\in Ver(f)$ and its normal part $\delta f^{\perp}\in Nor(f)$. In other words one as 
$$
T\mathcal{F} = Ver \oplus Nor
$$ 
as a direct sum of smooth fiber bundles over $\mathcal{F}$. Moreover this decomposition is preserved by the action of the re-parametrization group $\Gamma$, i.e. $\left(\delta f\circ\gamma\right)^{T} = 
 \delta f^{T} \circ\gamma$ and $\left( \delta f\circ\gamma\right)^{\perp} =  \delta f^{\perp}\circ\gamma$. 
\end{proposition}

\IEEEproof
The uniqueness of the decomposition into tangential and normal direction comes from the uniqueness of the decomposition of a vector in $\mathbb{R}^3$ into a tangent vector and normal vector to the 
surface. The smoothness of the decomposition is a consequence of the smoothness of the tangent and normal bundles.
To see that $\Gamma$ preserves the normal bundle, note that if $\gamma\in\Gamma$ reads $\gamma = \left(\gamma_1(u, v), \gamma_2(u,v)\right)$ in a chart, then $(f\circ\gamma)_u = f_u\circ\gamma~\frac{\partial \gamma_1}{\partial u} + f_v\circ\gamma~\frac{\partial \gamma_2}{\partial u}$ and $(f\circ\gamma)_v = f_u\circ\gamma~\frac{\partial \gamma_1}{\partial v} + f_v\circ\gamma~\frac{\partial \gamma_2}{\partial v}$, therefore 
$$
\left(f\circ\gamma\right)_u\times\left(f\circ\gamma\right)_v = f_u\circ\gamma\times f_v\circ\gamma \left(\frac{\partial \gamma_1}{\partial u} \frac{\partial \gamma_2}{\partial v} - \frac{\partial \gamma_2}{\partial u}
\frac{\partial \gamma_1}{\partial v} \right).
$$ 
It follows that the unit normal vector field to the parametrized surface $f\circ\gamma$ reads 
$$
\frac{\left(f\circ\gamma\right)_u\times\left(f\circ\gamma\right)_v}{\|\left(f\circ\gamma\right)_u\times\left(f\circ\gamma\right)_v\|} $$
$$
= \frac{ f_u\circ\gamma\times f_v\circ\gamma}{\| f_u\circ\gamma\times f_v\circ\gamma\|}\cdot
\frac{\left(\frac{\partial \gamma_1}{\partial u} \frac{\partial \gamma_2}{\partial v} - \frac{\partial \gamma_2}{\partial u}
\frac{\partial \gamma_1}{\partial v} \right)
}{|\frac{\partial \gamma_1}{\partial u} \frac{\partial \gamma_2}{\partial v} - \frac{\partial \gamma_2}{\partial u}
\frac{\partial \gamma_1}{\partial v} 
|} = n\circ\gamma,
$$
where in the last equality we have used that $\gamma$ preserves the orientation of $\mathbb{S}^2$.
Therefore $ \delta f\circ\gamma = (\delta f^{T} + \delta f^{\perp})\circ\gamma =  \delta f^{T} \circ\gamma + \delta f^{\perp}\circ\gamma$ with $\delta f^{T} \circ\gamma\in Ver(f\circ\gamma)$ and $ \delta f^{\perp}\circ\gamma\in Nor(f\circ\gamma$). The uniqueness of the decomposition then implies $ \delta f^{T} \circ\gamma = (\delta f \circ\gamma)^T$ and $\left( \delta f\circ\gamma\right)^{\perp} =  \delta f^{\perp}\circ\gamma$.

\section*{Proof of $\Gamma$-invariance of the elastic metric}
Now we will prove the fact that the elastic metric is invariant by  the group of orientation-preserving re-parametrizations $\Gamma =  \operatorname{Diff}^{+}(\mathbb{S}^2)$. This means that 
$$
\langle\!\langle h \circ \gamma, k\circ\gamma \rangle\!\rangle_{f\circ \gamma} = \langle\!\langle h, k\rangle\!\rangle_{f}.
$$
for $\gamma\in \Gamma$ and any tangent vectors $h, k$ at $f\in\mathcal{F}$.

Denote by $\tilde{f} := f\circ \gamma$. Set $(g,  {n}_{f}) := \Phi(f)$,  and $(\tilde{g},  {\tilde{n}}_{\tilde{f}})  := \Phi(\tilde{f})$. Define $\tilde{h}:= h \circ \gamma$ and $\tilde{k} = k\circ \gamma$.
Let us compute the volume form of the metric $\tilde{g}$. For any $s\in \mathbb{S}^2$, one has
$
\operatorname{Jac} \tilde{f}(s) = \operatorname{Jac} \left(f\circ \gamma\right)(s)  = (\operatorname{Jac} f)(\gamma(s))\cdot\operatorname{Jac}\gamma(s),
$
and
$
\tilde{g}(s) = \left(\operatorname{Jac}\gamma\right)^T g(\gamma(s)) \left(\operatorname{Jac}\gamma\right).
$
Therefore
\begin{align*}
\det\tilde{g}(s)  = \det (\operatorname{Jac}\gamma)^T \det g(\gamma(s)) \det \operatorname{Jac}\gamma \\= (\det \operatorname{Jac}\gamma)^2\det g(\gamma(s)),
\end{align*}
and
\begin{align*}
|\tilde{g}(s)|^{\frac{1}{2}} = \sqrt{\det\tilde{g}(s)}  = \sqrt{(\det \operatorname{Jac}\gamma)^2\det g(\gamma(s))} \\ = |\det \operatorname{Jac}\gamma|\, |g(\gamma(s))|^{\frac{1}{2}}.
\end{align*}
Let us now compute the first two terms of the elastic metric. Since
$$
\tilde{g}(s)^{-1} = \left(\operatorname{Jac}\gamma\right)^{-1} g(\gamma(s))^{-1} \left(\operatorname{Jac}\gamma^{T}\right)^{-1},
$$
and
$$
\delta\!\tilde{g}(s) = \left(\operatorname{Jac}\gamma\right)^{T} \delta\! g(\gamma(s)) \operatorname{Jac}\gamma,
$$
one has
\begin{align*}
\begin{array}{l}
\!\!\!\operatorname{Tr} \tilde{g}^{-1} \delta\!\tilde{g}(s) =\\
 \operatorname{Tr}\left[\operatorname{Jac}\gamma^{-1} g(\gamma(s))^{-1} \left(\operatorname{Jac}\gamma^{-1}\right)^{T}\left(\operatorname{Jac}\gamma\right)^T
 \delta\! g(\gamma(s)) \operatorname{Jac}\gamma\right] 
 \\
 = \operatorname{Tr}\left[\left(\operatorname{Jac}\gamma\right)^{-1} g(\gamma(s))^{-1} \delta\! g(\gamma(s)) \operatorname{Jac}\gamma\right] \\= \operatorname{Tr} g^{-1}\delta\! g(\gamma(s)).
 \end{array}
\end{align*}
Therefore, if one denotes by $(\delta\!\tilde{g}_1, \delta\!  {\tilde{n}}_1)$ (resp. $(\delta\!\tilde{g}_2, \delta\!  {\tilde{n}}_2)$) the infinitesimal variation of the pull-back metric $\tilde{g}$ and the normal vector field $n_{f\circ\gamma}$ induced by the tangent vector $\tilde{h}\in T_{f\circ\gamma}\mathcal{F}$ (resp. $\tilde{k}\in T_{f\circ\gamma}\mathcal{F}$), and $(\delta\!{g}_1, \delta\!  {{n}}_1)$ (resp. $(\delta\!{g}_2, \delta\!  {{n}}_2)$) the infinitesimal variation of the pull-back metric ${g}$ and the normal vector field $n_{f}$ induced by the tangent vector ${h}\in T_{f}\mathcal{F}$ (resp. ${k}\in T_{f}\mathcal{F}$), one has
$$
\begin{array}{l}
\operatorname{Tr} \tilde{g}^{-1} \delta\!\tilde{g}_1 \tilde{g}^{-1} \delta\!\tilde{g}_2(s)=\\
 \operatorname{Tr}\left[\left(\operatorname{Jac}\gamma\right)^{-1} g^{-1} \delta\! g_1 \left(\operatorname{Jac}\gamma\right) \left(\operatorname{Jac}\gamma\right)^{-1} 
 g^{-1} \delta\! g_2 \operatorname{Jac}\gamma\right]\\
 = \operatorname{Tr} g^{-1}\delta\! g_1 g^{-1}\delta\! g_2(\gamma(s)).
 \end{array}
 $$
 For the last term of the metric, since $\gamma$ acts by re-parametrization on the normal vector field, one has $ \delta\!  {\tilde{n}}_1(s) =  \delta\!  {n}_1(\gamma(s))$ and 
 $ \delta\!  {\tilde{n}}_2(s) =  \delta\!  {n}_1(\gamma(s))$. The invariance by re-parametrization of the elastic metrics then follows by a simple change of variables in the integral defining it.

\section*{Alignment of 3D-shapes}
 
\begin{figure}[!ht]
\centering
\vspace{.5cm}
\includegraphics[width=8cm]{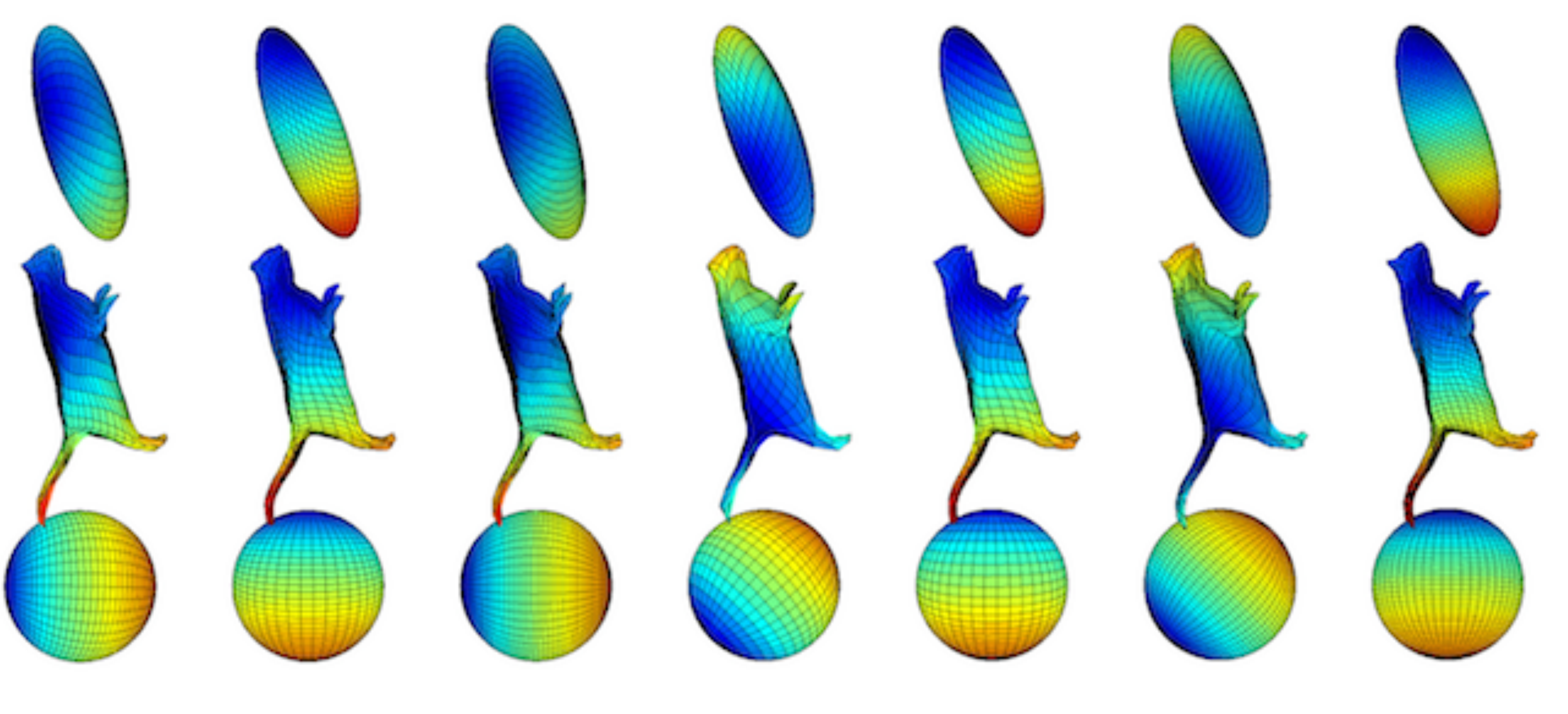}
\includegraphics[width=8cm]{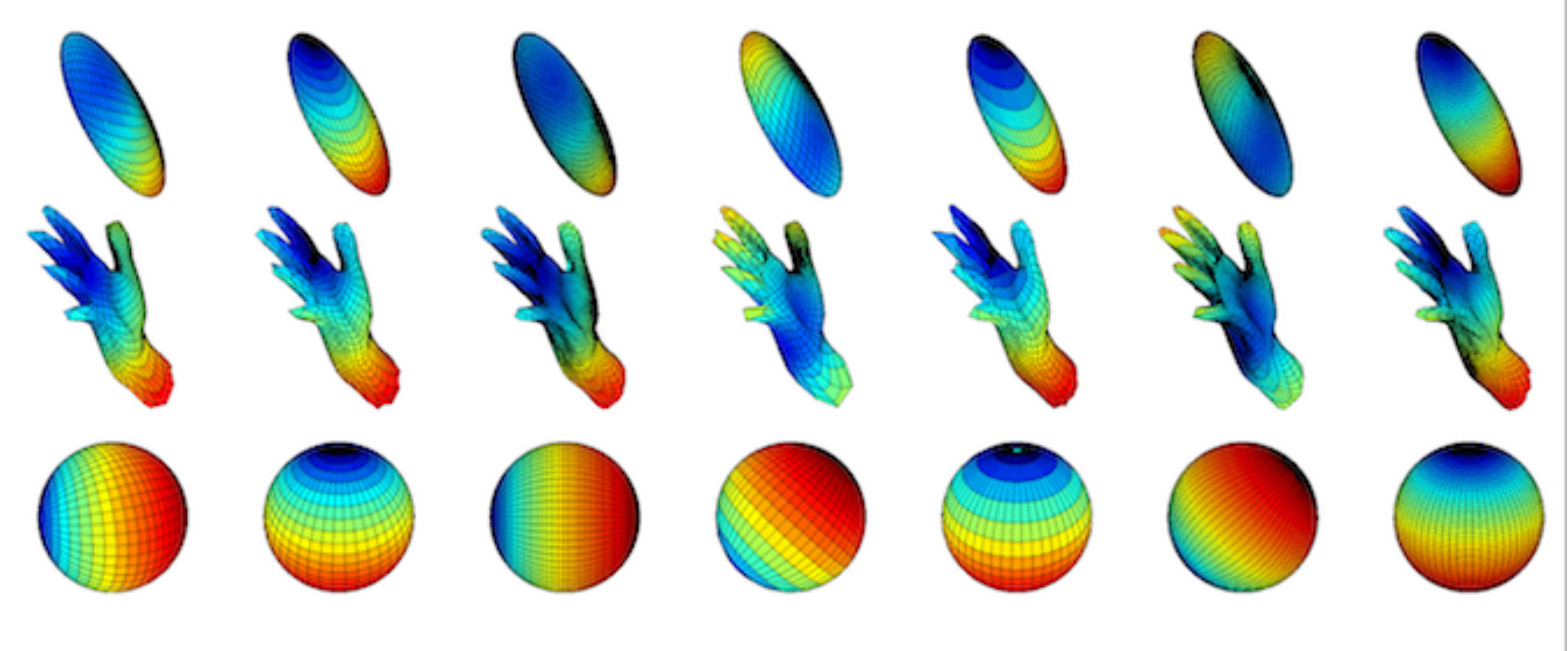}
\includegraphics[width=8cm]{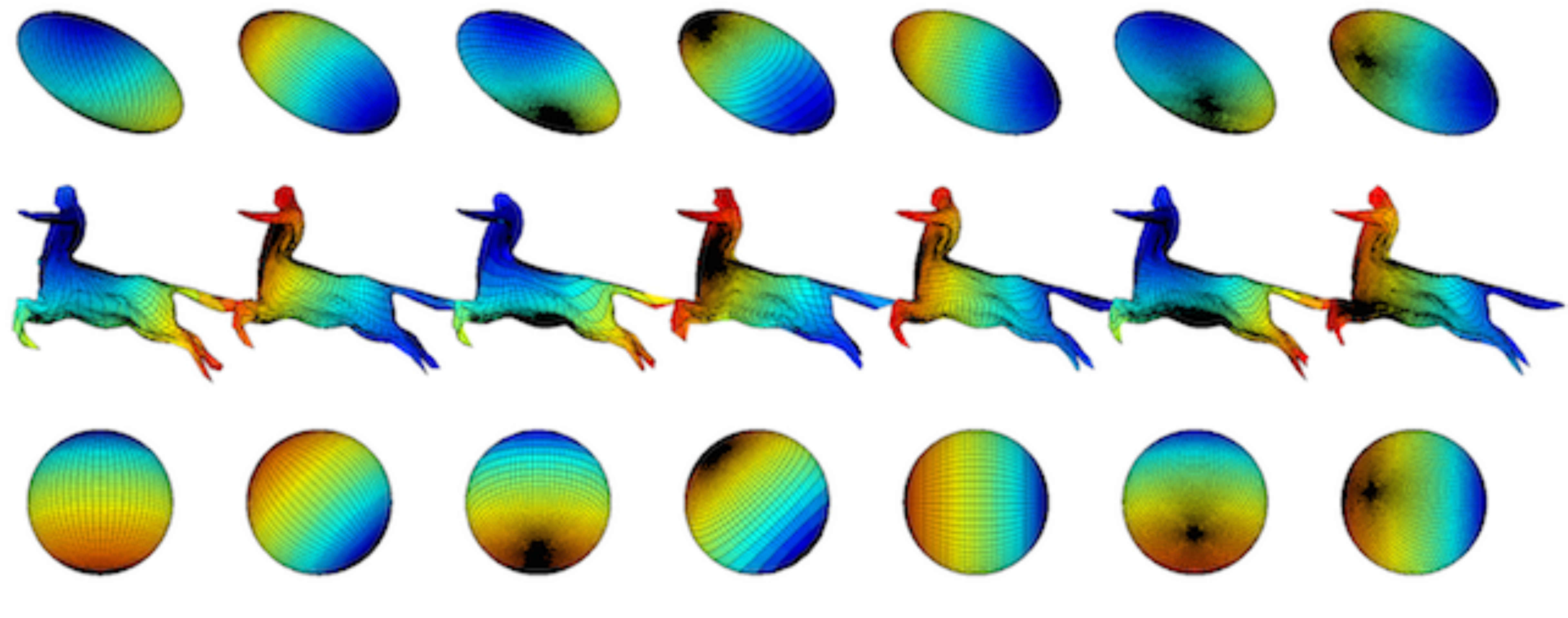}

\caption{Lower rows: different re-parametrizations of the sphere; Middle rows: corresponding re-parametrizations of a shape; Upper rows: corresponding approximating ellipsoids.}
\label{approx_ellipse}
\end{figure}

In this section we provide the details of the alignement program (Algorithm \ref{alg:Alignement}) described in section~4.1 of the paper.
The computation of the inscribed volume in a surface $f$ is made using Algorithm \ref{compute_volume}.
The center of mass of the inscribed volume in a surface $f$ is computed using Algorithm \ref{compute_center}.
The computation of the second moments is implemented using Algorithm \ref{compute_moments}.
More pictures illustrating the robustness of the approximating ellipsoid when the parametrization of the initial surface is changed are given in Fig.~\ref{approx_ellipse}. The dependance of the approximating ellipsoid with respect to a rotation of the surface is illustrated in Fig.\ref{rotation}.

\begin{figure}[!ht]
\centering
\includegraphics[width=8cm]
{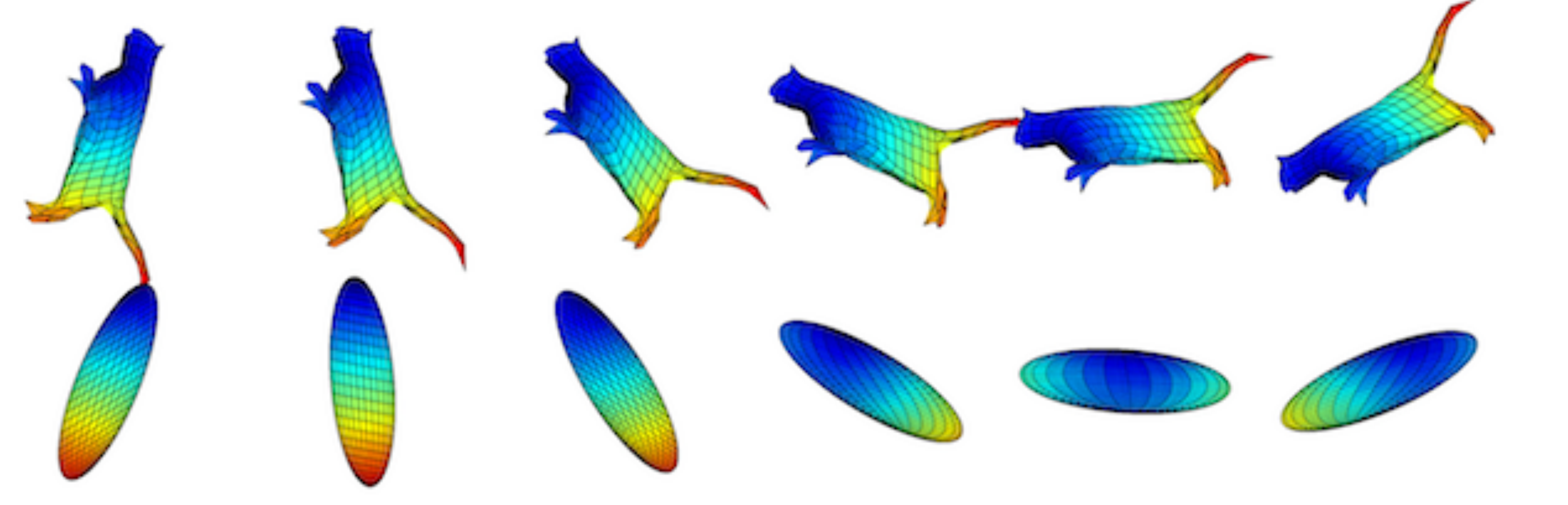}

\includegraphics[width=8cm]
{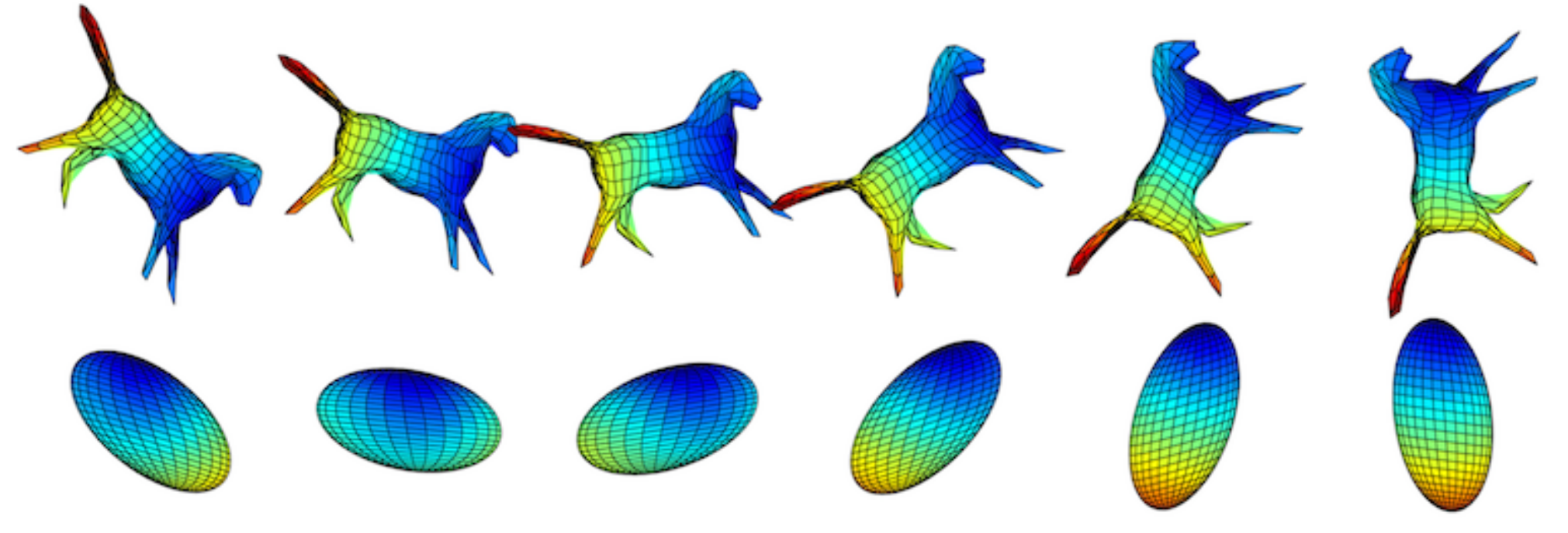}

\includegraphics[width=8cm]
{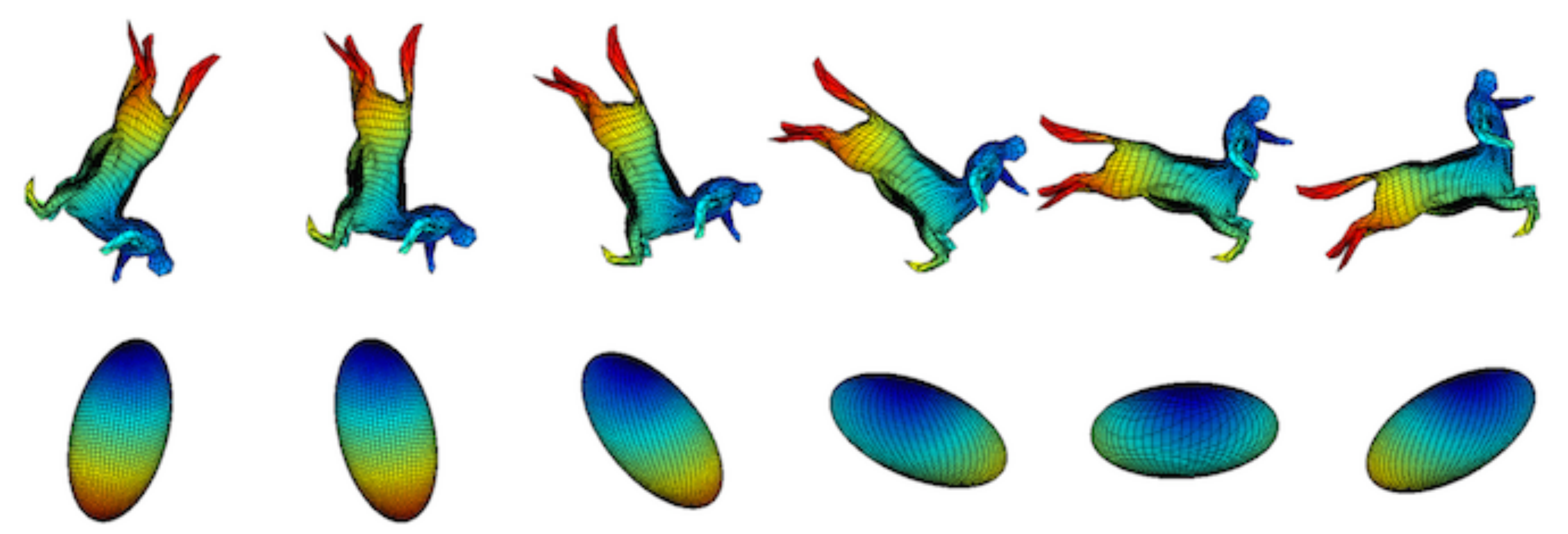}

\caption{Dependance of the approximating ellipsoid with respect to rotation of the shape.}
\label{rotation}
\end{figure}

\begin{algorithm}
\begin{footnotesize}
\KwIn{\begin{enumerate}
\item a grid of $n\times n$ points on the unit sphere, i.e. for each index $(i,j)\in [1,n]\times[1,n]$, a value of polar angle $\theta(i,j)$ and of azimuthal angle $\phi(i,j)$,
\item a parametrized surface $f_1$, i.e. for each index $(i,j)\in [1,n]\times[1,n]$, a point $f_1(i,j)$ in $\mathbb{R}^3$ corresponding to the image of the point on the sphere with spherical coordinates $(\theta(i,j), \phi(i,j))$ by the map $f_1$,
\item a parametrized surface $f_2$, i.e. for each index $(i,j)\in [1,n]\times[1,n]$, a point $f_2(i,j)$ in $\mathbb{R}^3$ corresponding to the image of the point on the sphere with spherical coordinates $(\theta(i,j), \phi(i,j))$ by the map $f_2$.
\end{enumerate}}
\KwOut{\begin{enumerate}
\item a centered and scaled surface $F_1$ having the same shape as $f_1$, with center of mass at the origin and inscribed volume $1$,
\item a centered, scaled and rotated surface $F_2$ having the same shape as $f_2$, with center of mass at the origin, inscribed volume $1$ and principal axes aligned with the principal axes of $F_1$.
\item For $k = 1,2$, an approximating ellipse $E_k$ of $F_k$.
\end{enumerate}}

\textbf{Algorithm:}
\begin{enumerate}
\item[\textbf{1-}] For $k=1,2$, use algorithm \ref{compute_volume} to compute  the volume $Vol_k$ inscribed in the surface $f_k$.

\item[\textbf{2-}] For $k=1,2$, $f_k \leftarrow f_k/\left(\textrm{Vol}_k\right)^{1/3}$.

\item[\textbf{3-}] For $k=1,2$, use algorithm \ref{compute_center} to compute the center of mass $\textrm{Center}_k$ of the inscribed volume in surface $f_k$.

\item[\textbf{4-}] For $k=1,2$, $f_k \leftarrow f_k - \textrm{Center}_k$. 

\item[\textbf{5-}] For $k=1,2$, use algorithm \ref{compute_moments} to compute the second moments $M_k$ of surface $f_k$.

\item[\textbf{6-}] For $k=1,2$, compute $[U_k,S_k, V_k] = \textrm{svd}(M_k)$.

\item[\textbf{7-}]
\textbf{Set} $F_1 = f_1$ and $F_2 = U_2\times U_1'\times f_2$.
\item[\textbf{8-}] For $k=1,2$, compute $$A_k = \left(\frac{4\pi}{15}\right)^{\frac{1}{5}}\det(M_k)^{-\frac{1}{10}}U_k\times\sqrt{S_k}\times U_k'.$$

\item[\textbf{9-}]
\textbf{Set} $E_k = A_k\times\textrm{sphere}$, $k=1,2$.
\end{enumerate}
\end{footnotesize}
\caption{\footnotesize Alignement of $3D$-shapes}
\label{alg:Alignement}
\end{algorithm}

\begin{algorithm}
\begin{footnotesize}
\KwIn{ $3D$-parametrized surface $f$ of size $a\times b\times 3$.}
\KwOut{\begin{enumerate}
\item
Inscribed volume $\textrm{Vol}$ in surface $f$ \item 
volume $\textrm{vol}{1}(i,j)$ of infinitesimal tetrahedron with vertices $0, f(i,j,:), f(i+1,j), f(i,j+1)$;
\item 
volume $\textrm{vol}{2}(i,j)$ of infinitesimal tetrahedron with vertices $0, f(i+1,j+1,:), f(i+1,j), f(i,j+1).$
\end{enumerate}}
\textbf{Algorithm:}
Initialize $\textrm{Vol} = 0$.

\For{$i\leftarrow 1$ \KwTo $size(f,1)$}{
\For{$j\leftarrow 1$ \KwTo $size(f,2)$}{
\begin{enumerate}
\item[\textbf{1-}]
\textbf{Set} \begin{itemize}
\item[]
$\textrm{edge}(1) = f(i+1, j,:)-f(i,j,:)$
\item[]
$\textrm{edge}(2) = f(i, j+1,:)-f(i,j,:)$
\item[]
$\textrm{edge}(3) = f(i, j+1,:)-f(i+1,j+1,:)$
\item[]
$\textrm{edge}(4) = f(i+1, j,:)-f(i+1,j+1,:)$
\end{itemize}
\item[\textbf{2-}] \textbf{Set} \begin{itemize}
\item[]$\textrm{vol}{1}(i,j) = \frac{1}{6}\textrm{Det}(\textrm{edge}(1), \textrm{edge}(2), -f(i,j,:))$
\item[]
\item[]
$\textrm{vol}{2}(i,j) = \frac{1}{6}\textrm{Det}(\textrm{edge}(3), \textrm{edge}(4), -f(i\!+\!1,j\!+\!1,:))$
\end{itemize}
\item[\textrm{3-}] $\textrm{Vol}\leftarrow \textrm{Vol} + \textrm{vol}(1) + \textrm{vol}(2).$
\end{enumerate}
}
}
\end{footnotesize}
\caption{\footnotesize Computation of inscribed volume}
\label{compute_volume}
\end{algorithm}

\begin{algorithm}
\begin{footnotesize}
\KwIn{ \begin{enumerate}
\item
$3D$-parametrized surface $f$ of size $a\times b\times 3$
\item
Inscribed volume $\textrm{Vol}$ in surface $f$ \item 
volume $\textrm{vol}{1}(i,j)$ of infinitesimal tetrahedron with vertices $0, f(i,j,:), f(i+1,j), f(i,j+1)$
\item 
volume $\textrm{vol}{2}(i,j)$ of infinitesimal tetrahedron with vertices $0, f(i+1,j+1,:), f(i+1,j), f(i,j+1).$
\end{enumerate}
}
\KwOut{Center of mass  of inscribed volume in surface $f$.}
\textbf{Algorithm:} Initialize $\textrm{Center} = (0,0,0).$
\For{$i\leftarrow 1$ \KwTo $size(f,1)$}{
\For{$j\leftarrow 1$ \KwTo $size(f,2)$}{
\begin{enumerate}
\item[\textbf{1-}]
$\textrm{m}{1} = \frac{1}{4}(f(i,j,:) + f(i+1,j,:)+f(i,j+1,:))$
\item[\textbf{2-}]
$\textrm{m}{2} = \frac{1}{4}(f(i+1,j+1,:) + f(i+1,j,:)+f(i,j+1,:))$
\item[\textbf{3-}]
$\textrm{Center} \leftarrow \textrm{Center} + 
\textrm{vol}{1}(i,j)\!\times\!\textrm{m}{1}\!
+ \!
\textrm{vol}{2}(i,j)\!\times\!\textrm{m}{2}
$
\end{enumerate}
}
}

\end{footnotesize}
\caption{\footnotesize Computation of center of mass}
\label{compute_center}
\end{algorithm}

\begin{algorithm}
\begin{footnotesize}
\KwIn{ \begin{enumerate}
\item
$3D$-parametrized surface $f$ of size $a\times b\times 3$
\item
Inscribed volume $\textrm{Vol}$ in surface $f$ \item 
volume $\textrm{vol}{1}(i,j)$ of infinitesimal tetrahedron with vertices $0, f(i,j,:), f(i+1,j), f(i,j+1)$
\item 
volume $\textrm{vol}{2}(i,j)$ of infinitesimal tetrahedron with vertices $0, f(i+1,j+1,:), f(i+1,j), f(i,j+1).$
\end{enumerate}
}
\KwOut{second moments of surface $f$ defined as the following integral over the inscribed volume
$$
M = \int \left(\begin{array}{ccc} x^2 &xy& xz\\ xy& y^2& yz\\ xz& xy& z^2\end{array}\right) \textrm{dvol}
$$
}
\textbf{Algorithm:} Initialize $M = \textrm{zeros}(3,3)$.

\For{$i\leftarrow 1$ \KwTo $size(f,1)$}{
\For{$j\leftarrow 1$ \KwTo $size(f,2)$}{
\For{$k\leftarrow 1$ \KwTo $3$}{
\For{$l\leftarrow 1$ \KwTo $3$}{
\tiny{
$$\textrm{s1} = (f(i,j,k)+f(i,j\!+\! 1,k))\!*\!(f(i,j,l)+f(i,j\!+\! 1,l)),$$

$$\textrm{s2} = (f(i,j,k)+f(i\!+\! 1,j,k))\!*\!(f(i,j,l)+f(i\!+\! 1,j,l)),$$

$$\textrm{s3} = (f(i\!+\!1,j,k)\!+\!f(i,j\!+\!1,k))\!*\!(f(i\!+\!1,j,l)\!+\!f(i,j\!+\!1,l)).$$

$$\textrm{m1} = \frac{1}{20}*(\textrm{s1}+\textrm{s2}+\textrm{s3}).$$


$$\textrm{s4} = (f(i\!+\!1,j,k)\!+\!f(i,j\!+\! 1,k))\!*\!(f(i\!+\! 1,j,l)\!+\!f(i,j\!+\! 1,l)).$$

$$\textrm{s5} = (f(i\!+\! 1,j\!+\! 1,k)\!+\!f(i\!+\! 1,j,k))\!*\!(f(i\!+\!  1,j\!+\! 1,l)+f(i\!+\! 1,j,l)),$$

$$\textrm{s6} = (f(i\!+\! 1,j\!+\! 1,k)\!+\! f(i,j\!+\! 1,k))\!*\!(f(i\!+\!1,j\!+\! 1,l)+f(i,j\!+\! 1,l)),
$$

$$\textrm{m2} = \frac{1}{20}
*(\textrm{s4}+\textrm{s5}+\textrm{s6}).$$

 $$M(k,l) \leftarrow M(k,l) + \textrm{vol1}(i,j).*\textrm{m1} +\textrm{vol2}(i,j).*\textrm{m2}.$$
 }
}
}
}
}

\end{footnotesize}
\caption{\footnotesize Computation of second moments}
\label{compute_moments}
\end{algorithm}






\newpage
\section*{Second order approximation of a surface}
 Algorithm \ref{algo:courbures} gives the second order approximation of a surface at a given point. It was used in section~4.2 in order to compute the principal curvatures of a surface.
\begin{algorithm}
\begin{footnotesize}
\KwIn{a surface passing through the origin, tangent to the $xy$-plane at the origin}
\KwOut{coefficients $A = (a_1,a_2,a_3,a_4,a_5,a_6)$ of the second order polynomial $P(x,y) = a_1x^2+a_2y^2+a_3xy+a_4x+a_5y+a_6$, which minimize the sum $\sum_i(z_i - P(x_i, y_i))^2$ over the points $(x_i, y_i,z_i)$ of the surface}
\textbf{Algorithm:}
Initialize $z_B = \textrm{zeros}(1,6)$, $B = \textrm{zeros}(1,6)$, $B2 = \textrm{zeros}(6,6)$.

\For{$i\leftarrow 1$ \KwTo $\textrm{number of points}$}
{\begin{enumerate}
\item[\textbf{1-}]\begin{itemize}
\item[]
$z_B(1) \leftarrow z_B(1) + z_ix_i^2;$
\item[]
$z_B(2) \leftarrow z_B(2) + z_iy_i^2;$
\item[]
$z_B(3) \leftarrow z_B(3) + z_ix_iy_i;$
\item[]
$z_B(4) \leftarrow z_B(4) + z_ix_i;$
\item[]
$z_B(5) \leftarrow z_B(5) + z_iy_i;$
\item[]
$z_B(6) \leftarrow z_B(6) + z_i;$
\end{itemize}
\item[\textbf{2-}]
\begin{itemize}
\item[]
$B(1) \leftarrow  x_i^2;$
\item[]
$B(2) \leftarrow  y_i^2;$
\item[]
$B(3) \leftarrow  x_iy_i;$
\item[]
$B(4) \leftarrow x_i$;
\item[]
$B(5) \leftarrow  y_i;$
\item[]
$B(6) \leftarrow  1;$
\end{itemize}
\item[\textbf{3-}]
$B2 \leftarrow  B2 + B'*B;$
\end{enumerate}
}
$A = \textrm{inv}(B2)*z_B'$.
\end{footnotesize}
\caption{
Computation of second order approximation of a surface at a given point.}
\label{algo:courbures}
\end{algorithm}

\newpage
\section*{Independance of the Energy function with respect to reparametrization}

From the theory it is clear that the energy function defined in the paper is independant of the way shapes are parametrized along a path. To provide numerical examples to illustrate this fact was however a difficult task. One has to mention here that the parametrization has to be changed \textit{smoothly} in order to provide a \textit{smooth} path in the pre-shape space. Therefore only parametrizations that are closed to the initial parametrization can be used in these experiments. 

\begin{figure*}[!ht]
 		\centering
		\includegraphics[width=16cm]
        {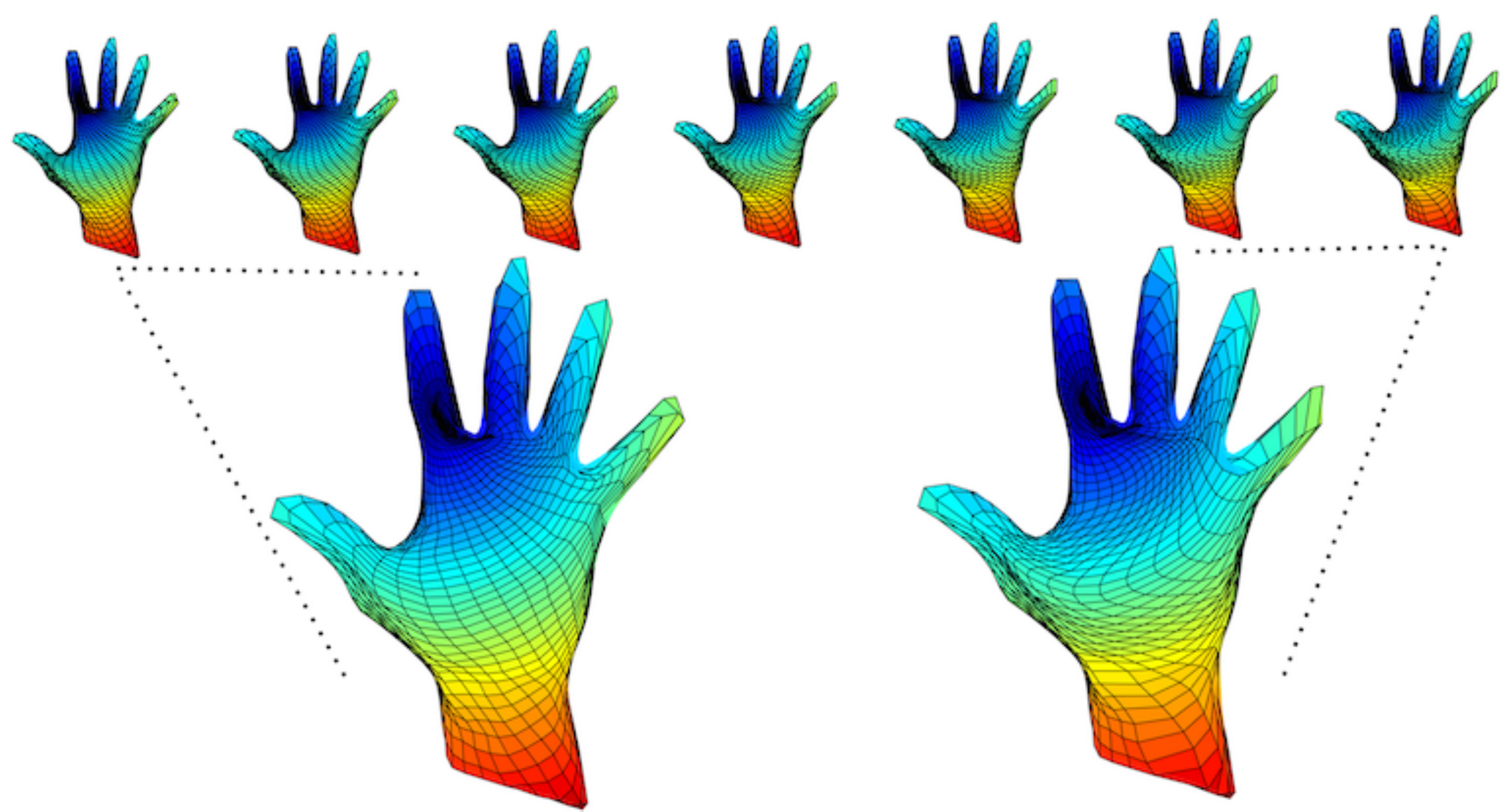}
        \includegraphics[width=15cm]
        {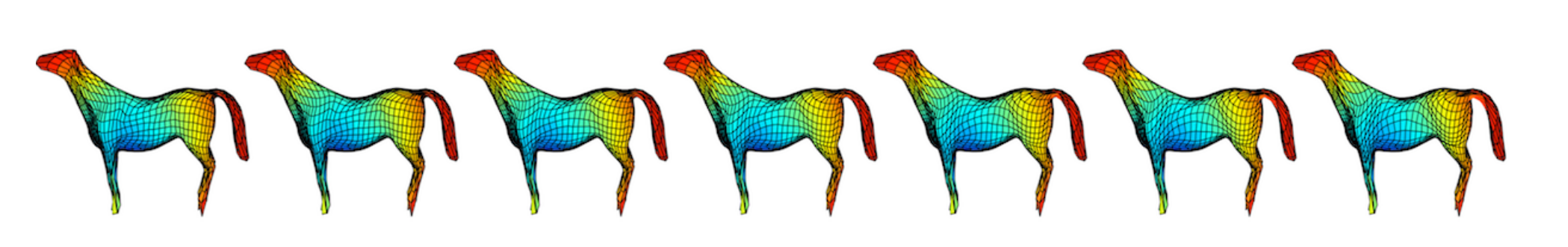}
        \includegraphics[width=15cm]
        {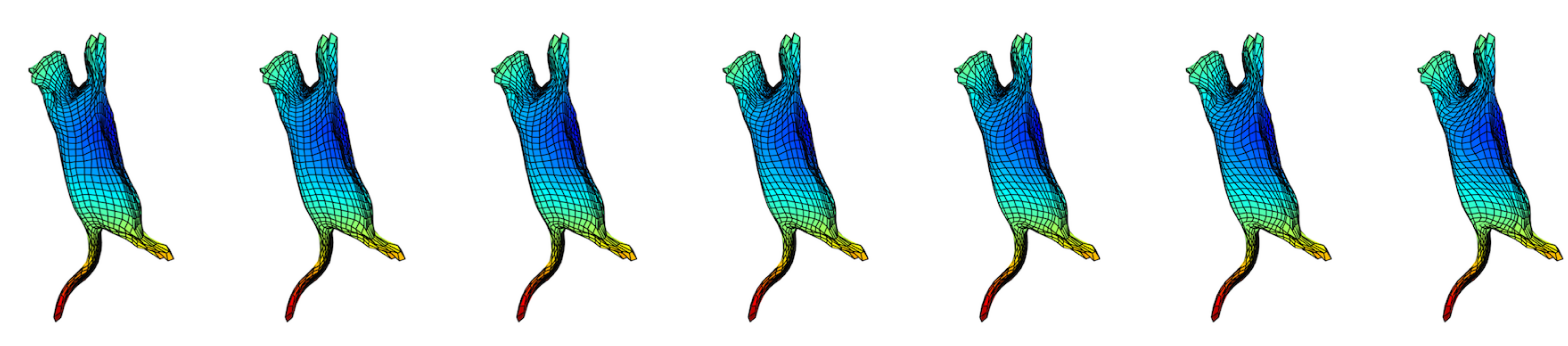}
        \includegraphics[width=15cm]
        {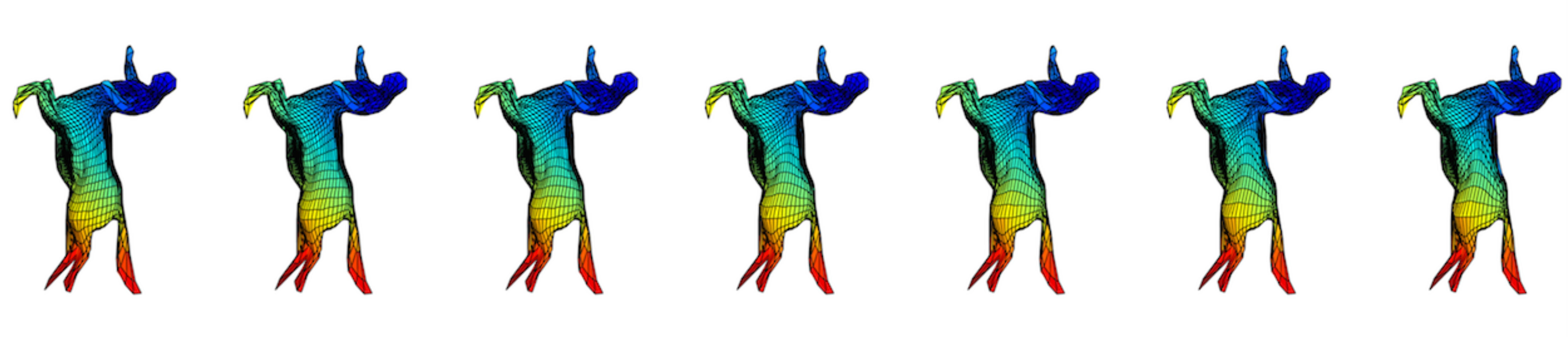}
 		\caption{Four Paths connecting the same shape but with a parametrization depending smoothly on time. The energy computed by our program is respectively $E_\Delta = 0$ for the path of hands, $E_\Delta = 0.1113$ for the path of horses, $E_\Delta = 0$ for the path of cats, and $E_\Delta = 0.0014$ for the path of Centaurs.}
 		\label{zero_energy}
 		\end{figure*}

In Figure \ref{zero_energy}, we give  examples of zero-energy paths, projecting to a point in shape space. The energy, as computed by our program, is closed to zero for each of them.

\begin{figure*}[!ht]
 		\centering
		\includegraphics[width=15cm]
{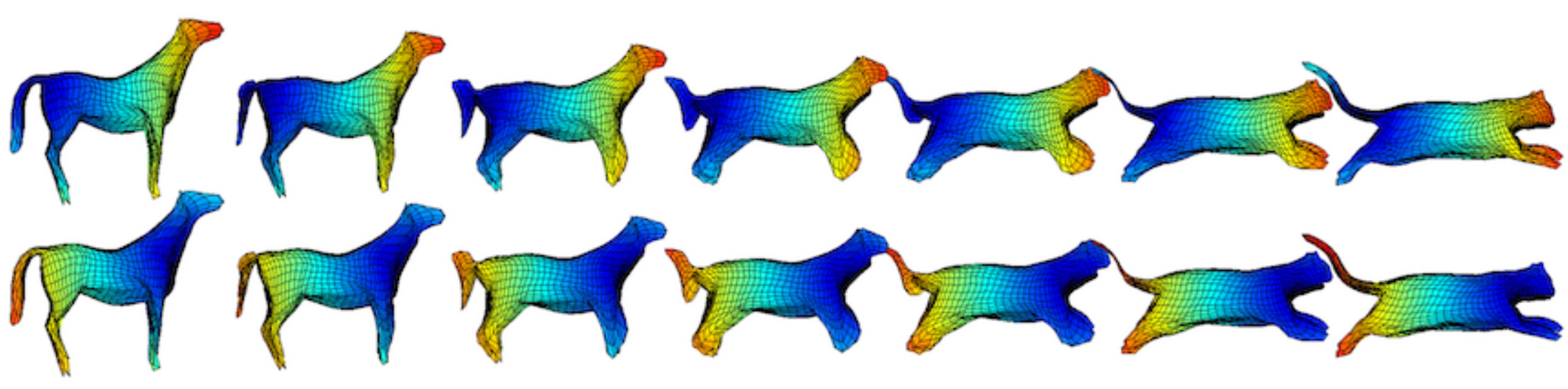}

\includegraphics[width=15cm]
{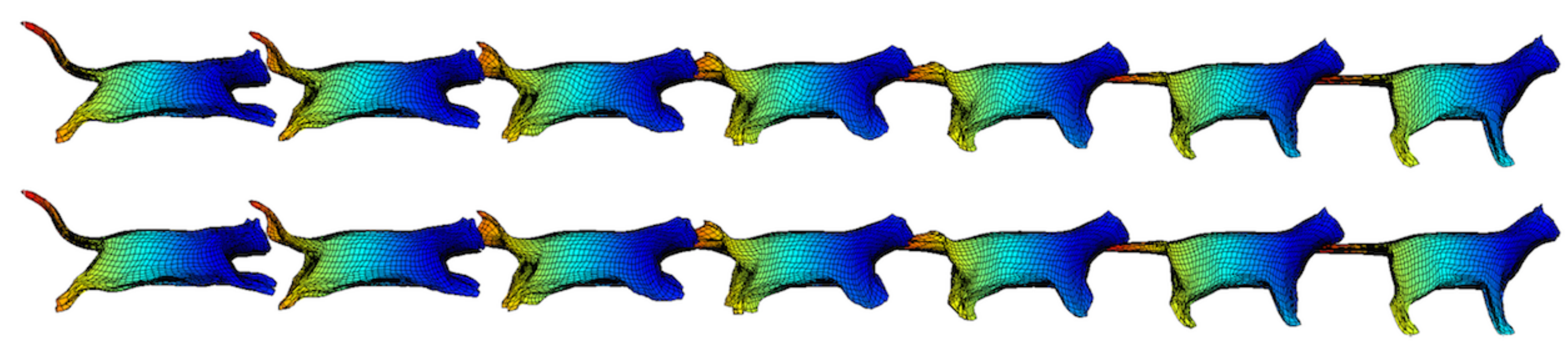}

\includegraphics[width=15cm]
{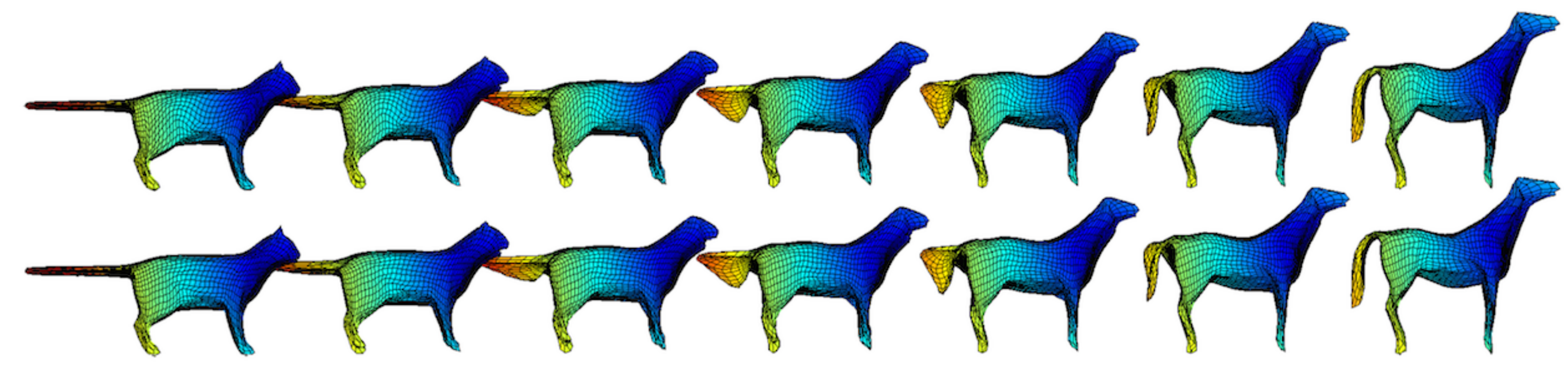}
\caption{Pairs of paths projecting to the same path in Shape space, but with different parametrizations. The energies of these paths, as computed by our program, are respectively (from the upper row to the lower row): $E_\Delta = 225.3565$, $E_\Delta = 225.3216$, $E_\Delta = 180.8444$, $E_\Delta = 176.8673$, $E_\Delta = 243.1812$ and $E_\Delta = 239.5410$.}
\label{same_energy}
 		\end{figure*}

In Figure \ref{same_energy}, we are interested in two different lifts of the same path in Shape space. 
The rows go by pairs: the upper two rows show a metamorphosis from a horse to a jumping cat, but with two different parametrizations. 
Theoretically the energy of the two upper paths should be the same. Numerically we obtain an energy $E_\Delta = 225.3565$ 
for the upper path, and $E_\Delta = 225.3216$ for the second one. For the third and forth paths, showing a metamorphosis 
from a jumping cat to a standing cat, the energy computed by our program is $E_\Delta = 180.8444$ and $E_\Delta = 176.8673$ respectively. For the fifth and sixth paths, from a standing cat to a standing horse, the computed energies
are $E_\Delta = 243.1812$ and $E_\Delta = 239.5410$ respectively. 
These energies were computed with the parameters $a=1$, $\lambda = 0.125$, $c=0$, $50^2$ numbers of points and using 
6-neighboordhoods for the computation of principal curvatures.

\newpage

\newpage

\begin{figure*}[!ht]
 		\centering
		\includegraphics[width=15cm]
        {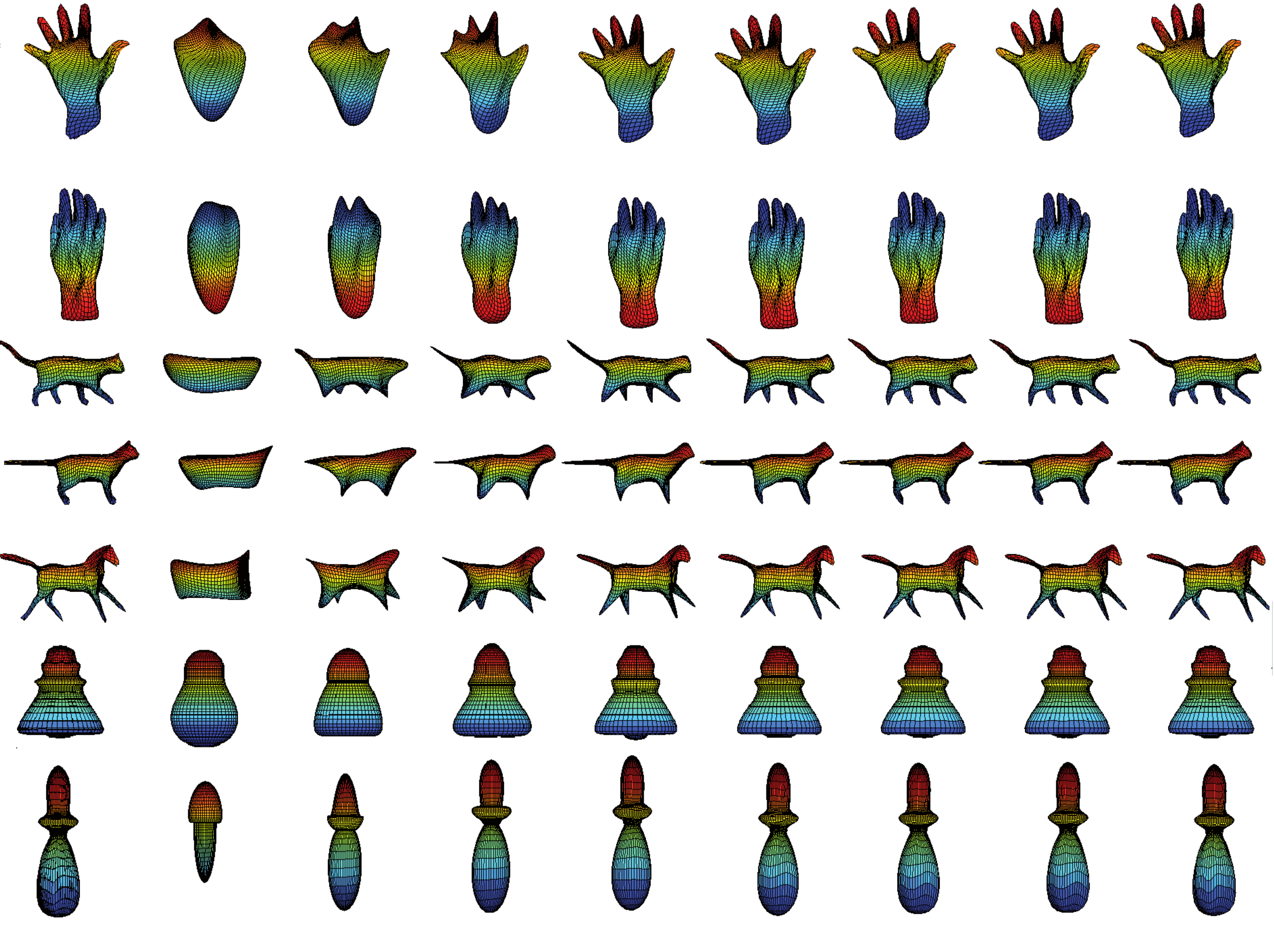}
 		\caption{Reconstruction of several surfaces with different degree of spherical harmonics. From left to right is depicted the initial surface and its approximation using spherical harmonics with maximal degree $l = 3$, $l = 5$, $l = 7$, $l=11$, $l=15$, $l=18$, $l=20$ and $l=28$ respectively.}
 		\label{reconstruction}
 		\end{figure*}

\begin{table*} [h]\label{tab:notation}
\footnotesize{
\begin{center} \caption{\bf List of symbols and their
definitions used in this paper.}
\begin{tabular}{||c|l||}
\hline \hline {\bf Symbol} & {\bf Definition} /{\bf Explanation} \\
\hline
$\mathbb{S}^2$ & the unit sphere in $\mathbb{R}^3$,\\
$(u, v)$ & local  coordinates on $\mathbb{S}^2$,\\
 $a \cdot b$ & the Euclidean inner dot product in
$\real^3$, \\
$a\times b$ & the cross product in $\mathbb{R}^3$.\\
\hline
$S$ & an observed surface modulo translation and rotation, or shape, \\
${\cal S}$ & the set of observed surfaces in $\real^3$ modulo translation and rotation or shape space,\\
$L(S(t))$& length of the path $t\mapsto S(t)$ of surfaces,\\
$dist(S_1, S_2)$ & infimum of the lengths of paths in $\mathcal{S}$ between $S_1$ and $S_2$,\\
$ T_S{\cal S}$ & the tangent space to the space of shapes $\mathcal{S}$ at a given shape $S$, \\ & It is identified with the space of all vector fields orthogonal to the surface $S$, \\
$\langle{X_1}, {X_2}\rangle$ & a Riemannian metric on the shape space ${\cal S}$, \\
$\phi_t(S_1,S_2)$ & a geodesic path in ${\cal S}$, from $S_1$ to
$S_2$, parameterized by $t \in [0,1]$, \\
& $\phi_0(S_1,S_2) = S_1,\ \phi_1(S_1,S_2) = S_2$. \\ \hline
$f$ or $F$ & a parametrized surface, i.e a smooth map $F~:\mathbb{S}^2\rightarrow \mathbb{R}^3$ \\& which is a homeomorphism onto its image and whose differential is injective,\\
$f_u$ and $f_v$ & derivatives of $f$ with respect to the coordinates $u$ and $v$ respectively,\\
${\cal F}$ & the set of parametrized surfaces in $\real^3$ or pre-shape space,\\
$L(F(t))$& length of the path $t\mapsto F(t)$ of parametrized surfaces,\\
$dist(F_1, F_2)$ & infimum of the lengths of paths in $\mathcal{F}$ between $F_1$ and $F_2$,\\
$ T_f{\cal F}$ & the tangent space to the pre-shape space $\mathcal{F}$ at a given parametrized shape $f$, \\ 
$\langle\!\langle{X_1}, {X_2}\rangle\!\rangle$ & a Riemannian metric on the pre-shape space ${\cal F}$, \\
$\Phi_t(F_1,F_2)$ & a geodesic path in ${\cal F}$, from $F_1$ to
$F_2$, parameterized by $t \in [0,1]$, \\
& $\phi_0(F_1,F_2) = F_1,\ \phi_1(F_1,F_2) = F_2$, \\
$B_{f}(X_1, X_2)$ &  a symmetric bilinear map on the tangent space $T_{f}\mathcal{F}$ which is non-negative \\&when apply to the same tangent vector,\\
$\textrm{Ker}B_f$ & Kernel of the symmetric bilinear map $B_f$ defined as\\ &
$\textrm{Ker}B_f = \{X_1\in T_f\mathcal{F}\textrm{~such~that~} B_f(X_1, X_2) = 0,~ \forall X_2\in T_f\mathcal{F}\}$,\\
$ (\!({X_1}, {X_2})\!)$ & a non-negative semi-definite inner product on $\mathcal{F}$, i.e. a symmetric bilinear map $B_f$\\& on each tangent space $T_{f}\mathcal{F}$
  which is non-negative when apply to the same tangent vector.\\
 \hline
$\textrm{SO}(3)$ & the group of rotations in $\mathbb{R}^3$,\\
$\textrm{SO}(3)\rtimes \mathbb{R}^3$ & the group of rotations and translations in $\mathbb{R}^3$,\\
$\textrm{Diff}(\mathbb{S}^2)$ & the group of all diffeomorphisms of  $\mathbb{S}^2$, i.e. reparametrization group,\\
$\Gamma = \textrm{Diff}^+(\mathbb{S}^2)$ & the group of orientation-preserving diffeomorphisms of  $\mathbb{S}^2$ \\&i.e. group of  reparametrizations preserving orientation,\\
$G =  \textrm{Diff}^+(\mathbb{S}^2)\times \textrm{SO}(3)\rtimes \mathbb{R}^3$ & shape-preserving group.\\
\hline
$\mathcal{F}/G$ & quotient space of the pre-shape space by the group of  reparametrizations preserving orientation, \\
& is identified with the shape space $\mathcal{S}$ by the application which maps \\&a parametrized surface $f~:\mathbb{S}^2\rightarrow \mathbb{R}^{3}$ to its image,\\ 
$[f]$ & an element of the quotient space $\mathcal{F}/G$
which is 
 the orbit of a parametrized shape $F$ \\& under the reparametrization group $\Gamma=\textrm{Diff}^+(\mathbb{S}^2)$ and the group $\textrm{SO}(3)\rtimes \mathbb{R}^3$\\&
 $[f] = \{(Of + v)\circ \gamma \textrm{~for~}O \in \textrm{SO}(3), v\in\mathbb{R}^3, \gamma \in \Gamma\}$ \\

\hline \hline
\end{tabular}
\end{center}
}
\end{table*}

\end{document}